\documentclass[hidelinks,onefignum,onetabnum]{siamart220329}

\usepackage{amsfonts,upquote}
\usepackage{graphicx,verbatim}

\headers{Quadrature formulas from rational approximation}{Horning and Trefethen}

\title{{Quadrature formulas from rational approximations}%
\thanks{Submitted to the editors DATE.}}

\author{Andrew Horning\thanks{Dept.\ of Mathematical Sciences,
Rensselaer Polytechnic Institute, Troy, NY 12180, USA
(\email{hornia3@rpi.edu})} \and Lloyd N. Trefethen\thanks{School of Engineering and Applied Sciences,
Harvard University, Cambridge, MA 02138, USA}
(\email{trefethen@seas.harvard.edu})}

\def\complex{{\mathbb{C}}}
\def\Re{\hbox{\rm \kern .8pt Re\kern .6pt}}
\def\Im{\hbox{\rm \kern .8pt Im\kern .6pt}}
\def\half{{1\over 2}}

\def\zk{z_k^{}}
\def\ck{c_k^{}}
\def\rn{r_{\kern -.5pt n}^{}}
\def\In{I_n^{}}
\def\ti{{1\over 2\pi i}}
\def\riS{r^{-1}\kern -1pt S}
\def\C{{\kern .7pt\cal C}}
\def\negreal{(-\infty,0\kern .3pt]}
\def\pput(#1,#2)#3{\noindent\smash{\raise#2pt\hbox to 0pt
   {\kern #1pt #3\hss}}\ignorespaces}

\begin{document}

\maketitle

\begin{abstract}
It is shown that quadrature formulas in many different applications
can be derived from rational approximation of the Cauchy transform of
a weight function.
Since rational approximation
is now a routine technology, this provides an easy new method to derive
all kinds of quadrature formulas as well as fundamental insight into
the mathematics of quadrature.
Intervals or curves of quadrature nodes correspond
to near-optimal branch cuts of the Cauchy transform.
\end{abstract}



\section{\label{intro}Introduction}

Quadrature formulas are used across scientific computing,
and in numerical linear algebra, they have become an
important tool for reducing large-scale eigenvalue and
matrix function problems to small sets of matrix solves
\cite{beyn,beg,bst,col,ct,elguide,goedecker,gptv,gt,nicks,hg,polizzi,ss}.
Traditionally they are derived and analyzed on a case-by-case
basis.  For one example, Gauss quadrature is connected with the
theory of orthogonal polynomials~\cite{gautschi}.  For another,
the exponentially convergent periodic trapezoidal rule is
analyzed by discrete Fourier analysis or estimation of contour
integrals~\cite{trap}.  More complicated quadrature problems are
reduced to known ones such as these with the aid of conformal
maps~\cite{gt,haletee,ht,bromwich}, as in the ``three Nicks''
paper by Nick Higham, Nick Hale, and the second author~\cite{nicks}
(see Figure~\ref{nicksfig}).

In this article we show how quadrature formulas can be
generated with very little effort by rational approximation,
often approximation of functions with a ``$\hbox{sign}(z)$'' or
``Zolotarev'' flavor~\cite{it,zolo}.  The quadrature nodes are
the poles of the rational function, delineating an approximate
branch cut on a real interval or a complex contour, and the
quadrature weights are the residues.  This gives new insight into
the mathematics of quadrature formulas.  Until recently, there
would have been few computational implications, because computation
of rational approximations was a difficult problem, but with the
appearance of the AAA (adaptive Antoulas--Anderson) algorithm in
2018~\cite{aaa,acta}, together with the AAA-Lawson enhancement to
upgrade near-best to best~\cite{lawson}, it has become routine.
The Chebfun code {\tt aaa.m}, which runs in MATLAB or
Octave~\cite{chebfun}, is used for the examples of this paper,
and AAA implementations are also available in Julia~\cite{driscoll}
and SciPy \cite{virtanen}.  More recently another
approximation algorithm due to Salazar Celis, based on Thiele
continued fractions rather than barycentric representations, has
also shown great promise~\cite{driscoll,salazar}.

Section 2 lays out the mathematics of the connection between rational
approximation and quadrature, and 
sections 3--7 and 9--10 apply these ideas to
seven established quadrature problems, with section 8 summarizing
certain necessary approximation matters related to Zolotarev.
In each case, we first apply the rational approximation idea to a
standard configuration and then we show how it can be modified to
derive a quadrature formula for a variant problem.  Section~11 closes
the paper with a discussion.

We must say a word about our personal ``journey.''  Once the schema of
section 2 became clear, it simply amazed us how, for one problem after
another, we were able to use it to reproduce, with a fraction of a
second of computing time, results won previously by careful analysis:
\medskip

{\obeylines\leftskip 1.5cm
     Figure 3: see Fig.\ 4.1 of \cite{ht},
     Figure 4: see Fig.\ 1 of \cite{tee},
     Figure 5: see Fig.\ 2.3 of \cite{haletee},
     Figure 9: see Fig.\ 3 of \cite{deano09},
     Figure 10: see Fig.\ 4.3 of \cite{talbot},
     Figure 11: see Fig.\ 1 of \cite{bromwich},
     Figure 13: see Fig.\ 2 of \cite{akt},
     Figure 16: see Fig.\ 2 of \cite{nicks}.
\par}
\medskip

\section{\label{math} Mathematics of the approximation-quadrature connection}
\begin{figure}
\begin{center}
\includegraphics[trim=20 180 20 17, clip, scale=.80]{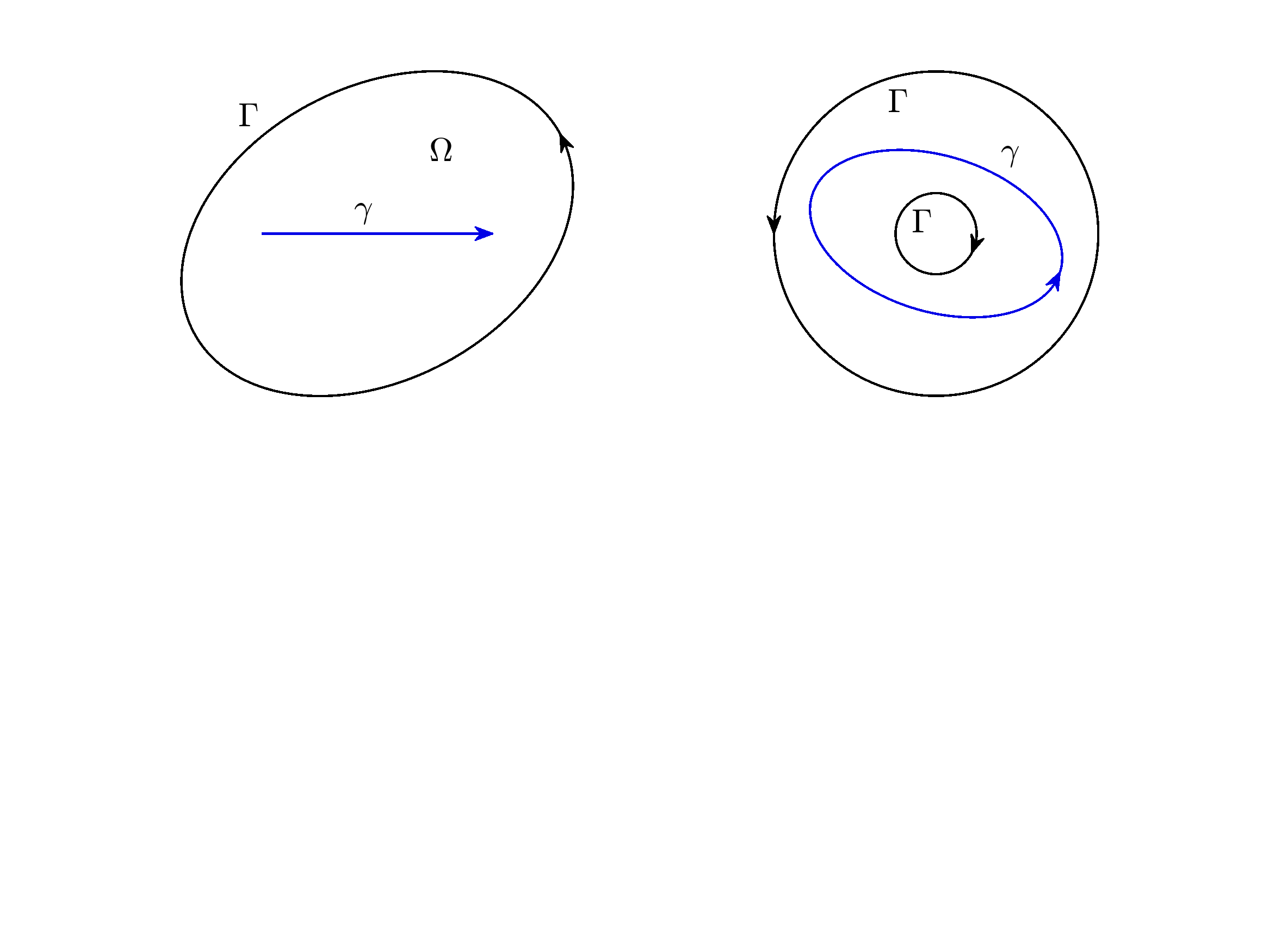}
\end{center}
\vspace*{-5pt}
\caption{\label{schematic} On the left,
sketch of an integration arc\/ $\gamma$ and a
surrounding contour\/ $\Gamma$ in the complex plane.
On the right, the analogous configuration for a closed integration
contour $\gamma$.}
\end{figure}
As on the left in Figure 1, 
let $\gamma$ be a Jordan arc in the complex plane $\complex$, that is,
the image of $[-1,1]$ under a real or complex homeomorphism.
The simplest example is $[-1,1]$ itself.
Let $w(z)$ be an integrable weight function defined on
$\gamma$, and suppose we are interested in approximating
the integral
\begin{equation}
I = \int_\gamma f(z) \kern 1pt w(z)\kern 1pt  dz .
\label{defn}
\end{equation}
Assume $f$ is analytic in the closure of a Jordan region $\Omega$ containing $\gamma$,
bounded by a Jordan curve $\Gamma$.
By the Cauchy integral formula, we have
\begin{equation}
I = \int_\gamma \ti \left[ \kern 1.5pt \int_\Gamma {f(s) \kern 1pt
ds\over s-z}\right] \kern -1pt w(z)\kern 1pt dz ,
\label{twoint1}
\end{equation}
hence by exchanging orders of integration,
\begin{equation}
I = \ti \int_\Gamma f(s) \kern -1pt \left[ \kern 2pt
\int_\gamma {w(z)\kern 1pt dz\over s-z}\right] ds .
\label{twoint}
\end{equation}
Let $\C(s)$ denote the function in square brackets,
the {\em Cauchy transform\/} of $w$,
\begin{equation}
\C(s) = \int_\gamma {w(z)\kern 1pt dz\over s-z},
\label{cauchytrans}
\end{equation}
which is analytic throughout 
$\complex\cup \{\infty\}\backslash \gamma$.
(This is also called the {\em Stieltjes transform,} and often
there is a factor $1/(2\pi i)$ in front.)
Then (\ref{twoint}) becomes
\begin{equation}
I = \ti \int_\Gamma f(s)\kern 1pt  \C(s)\kern 1pt  ds ,
\label{twoint2}
\end{equation}
and thus we have converted a real or complex integral over the arc $\gamma$
to a complex contour integral over the enclosing contour $\Gamma$.
Now suppose that $\rn$ is a rational function of degree $n$
that approximates $\C$ on $\Gamma$.
(A rational function of degree $n$ is one that can be written as
the ratio of two polynomials of degree at most~$n$.)
Then (\ref{twoint2}) suggests approximating $I$ by
\begin{equation}
\In = \ti \int_\Gamma f(s) \kern 1pt \rn(s) \kern 1pt ds .
\label{approx}
\end{equation}
Assume the poles of $\rn$ are $n$ distinct finite numbers
$z_1^{},\dots,z_n^{}$, so $\rn$ can be written
\begin{equation}
\rn(s) = C + \sum_{k=1}^n {\ck\over s-\zk},
\end{equation}
where $C = \rn(\infty)$ and $\ck$ is the residue of $\rn$ at $\zk$.
Then, assuming all the poles are enclosed
by $\Gamma$, residue calculus converts (\ref{approx}) to 
\begin{equation}
\In = \sum_{k=1}^n \ck \kern 1pt f(\zk) \kern 1pt.
\label{approx2}
\end{equation}
This is our quadrature approximation to the integral $I$.  If
\begin{equation}
\|\kern .5pt \C-\rn\|_\Gamma^{} \le \varepsilon,
\label{approxeps}
\end{equation}
where $\|\cdot\|_\Gamma^{}$ is the $\infty$-norm on $\Gamma$, then
by (\ref{twoint2}) and (\ref{approx}),
\begin{equation}
| I - \In| \le {\varepsilon \over 2\pi} \kern 1pt|\Gamma| \kern 1.5pt \|f\|_\Gamma^{},
\label{bound}
\end{equation}
where $|\Gamma|$ is the arc length of $\Gamma$.
If $|\Gamma|$ is large or infinite, it is often possible
to remove this factor while replacing
$\|\C-\rn\|_\Gamma^{}$ or $\|f\|_\Gamma^{}$ by the
corresponding $1$-norm.

The same reasoning just outlined also applies if $\gamma$ is a closed
Jordan contour, that is, a complex homeomorphism of the unit circle.
All that changes now is that for $\Gamma$ to ``enclose'' $\gamma$,
it should consist of two Jordan contours, one inside $\gamma$
and the other outside, as sketched on the right in Figure $1$.

\section{Application 1: Gauss--Legendre quadrature}

The most basic quadrature problem is integration of
a function $f$ over the interval $\gamma = [-1,1]$ with weight function $w(z)=1$:
\begin{equation}
I = \int_{-1}^1 f(z) \kern 1pt dz .
\label{gauss-leg}
\end{equation}
This is the original setting of Gauss quadrature, which is also
called Gauss--Legendre quadrature to distinguish it from cases with
non-constant weight functions (sections 5 and 6).  Gauss derived his method by
rational approximation of the Cauchy transform (\ref{cauchytrans})
of $w(z)=1$ on $[-1,1]$, which is $\C(s) = \log((s+1)/(s-1))$
\cite{gauss,gautschisurvey}.  The only difference from the mathematics of the last
section is that he used a Pad\'e approximation at $s=\infty$ rather
than an approximation over a contour $\Gamma$.  This was before
Cauchy integrals and long before Pad\'e, so Gauss's formulation,
emphasizing continued fractions and hypergeometric series, is quite different from ours.

To illustrate the derivation of essentially the same
quadrature formulas via rational approximation on a contour $\Gamma$,
consider the test integrand
$f(z) = 1/(1+20 z^2)$, which is analytic in the complex $z$-plane except at
$z = \pm i/\sqrt{20}$.  Motivated by standard results in
approximation and quadrature, let us design a quadrature rule for
functions analytic inside the 
Bernstein ellipse $\Gamma$ (ellipse with foci $\pm 1$)
passing through these two points \cite{atap}.  In the standard notation, this is the
ellipse with
parameter $\rho = 1/\kern -1pt\sqrt{20} + \sqrt{21/20} \approx 1.248$ (sum of
semiminor and semimajor axis lengths).  
For AAA approximation, the degree~$n$ can be prescribed, or a tolerance for
(\ref{approxeps}) can be specified so that $n$ is determined
adaptively.
Here we compute
an approximation $\rn(s) \approx \C(s)$ of degree $n=20$ 
in about 0.05 s on a laptop with the code
\medskip

{\small
\begin{verbatim}
       rho = 1/sqrt(20) + sqrt(21/20);
       c = rho*exp(2i*pi*(1:200)'/200);
       Z = (c+1./c)/2;
       F = log((Z+1)./(Z-1));
       [r,pol,res] = aaa(F,Z,'degree',20,'sign',1);
\end{verbatim}
\par}
\medskip

\begin{figure}
\begin{center}
\includegraphics[trim=0 140 20 7, clip, scale=.88]{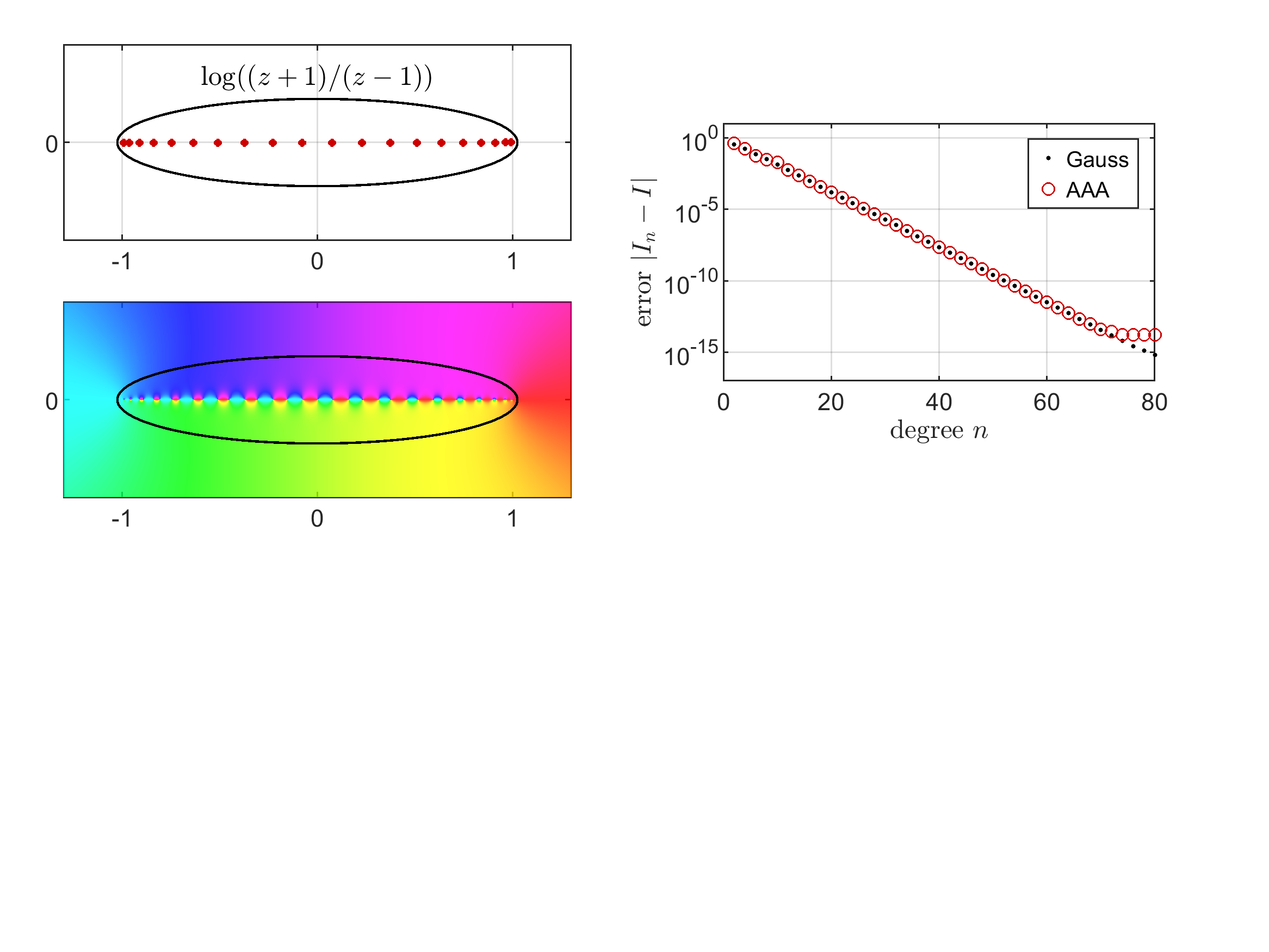}
\vspace*{-5pt}
\caption{\label{gaussfig}Quadrature by AAA rational approximation of\/ 
$C(s) = \log((z+1)/(z-1)))$ on the Bernstein ellipse passing through the
singularities at $\pm i/\sqrt{20}$
of the test integrand $f(z) = 1/(1+20z^2)$ on $[-1,1]$.  The 
upper-left image shows AAA poles for degree $n=20$, 
close to the Gauss quadrature nodes of the same degree, and the
lower-left image is a phase portrait of the AAA approximation $r_{20}^{}(z)$
(red for positive real, cyan for negative real).
The right image shows convergence for this integrand
as a function of degree $n$.}
\end{center}
\end{figure}

\noindent 
(The \verb|'sign'| flag is discussed in section 8.
As for the use of 200 sample points, the discretization details in
these computations are usually not important, and one can bypass them
by using the alternative continuum AAA algorithm \cite{driscoll,continuum}.)
The three outputs from {\tt aaa} are {\tt r}, a function
handle for the rational approximation, and {\tt pol} and {\tt res}, $n$-vectors
of its poles and corresponding residues.
Figure~\ref{gaussfig} shows that the poles look like
Gauss quadrature points, and one can check that the weights are
also close to Gauss weights.  Following (\ref{approx2}),
a quadrature result can be obtained with
\medskip

\begin{figure}
\begin{center}
\includegraphics[trim=0 140 20 7, clip, scale=.92]{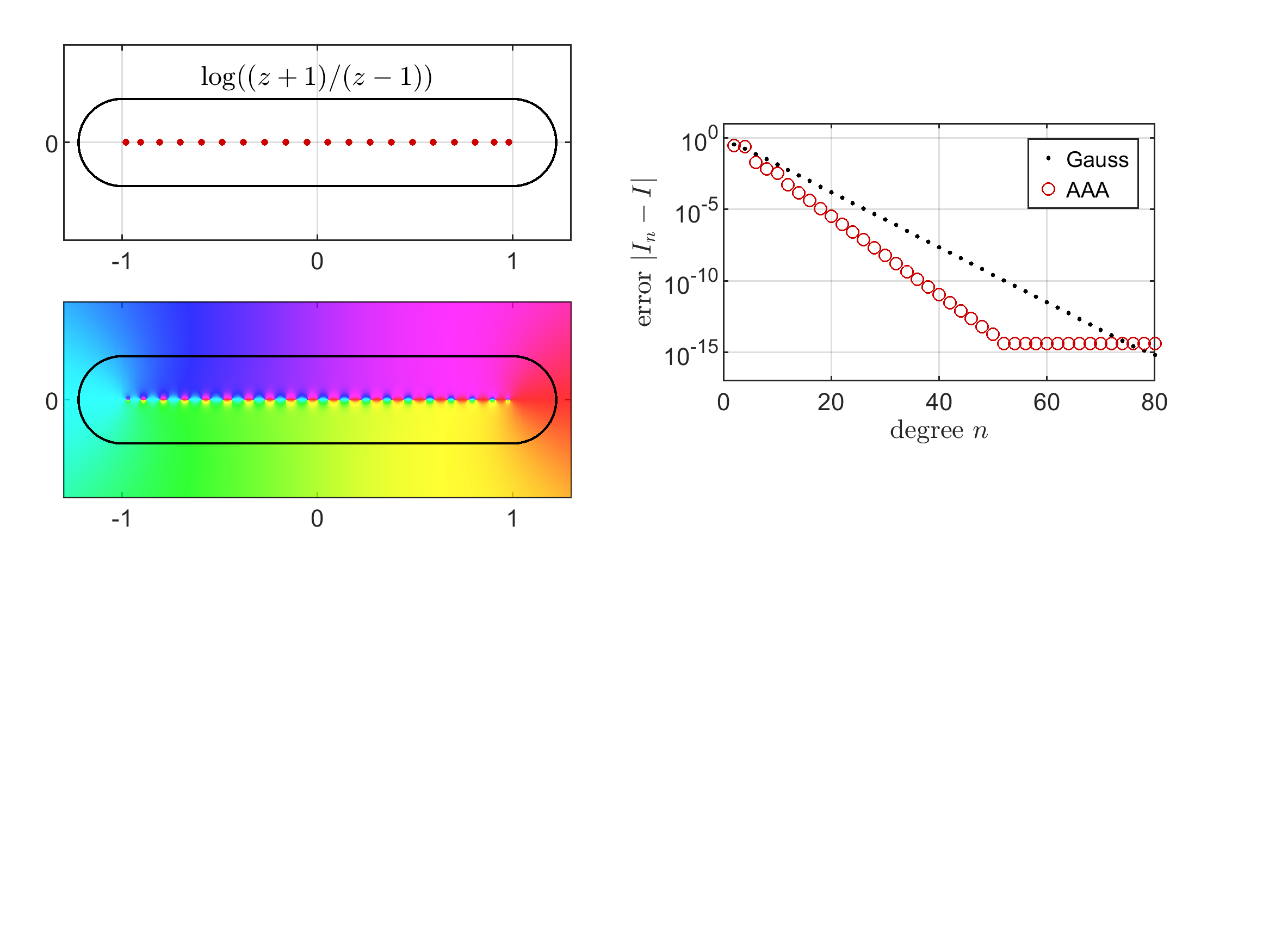}
\vspace*{-5pt}
\caption{\label{gaussvfig}Variant of Figure $\ref{gaussfig}$ in which
$\log((z+1)/(z-1))$ is approximated on an\/ $\varepsilon$-neighborhood
of\/ $[-1,1]$ rather than a Bernstein ellipse.  For the same test integrand as
before, the convergence is now faster by a factor of about $\pi/2$,
and the quadrature nodes are more evenly spaced.}
\end{center}
\end{figure}

{\small
\begin{verbatim}
       f = @(z) 1./(1+20*z.^2);
       In = res.'*f(pol)
\end{verbatim}
\par}
\medskip

\noindent 
and this gives ${\tt In} \approx
0.603943 - 0.00000002\kern .5pt i$, matching
the exact value $I\approx 0.604100$  
with an error of $1.6\cdot 10^{-4}$.

In this and other experiments of this paper, the nodes and
weights and computed integrals are all
slightly complex even when one would expect them to be real.
This is because the Chebfun implementation of AAA does
not enforce real symmetry, although this can be fixed, as is
done in some other implementations such as the {\tt rational} code
in the MathWorks RF Toolbox \cite{RF}.  With such an adjustment,
real symmetric problems give exactly real results.

The phase portrait in Figure \ref{gaussfig} shows a clear approximation to the branch
cut of $\log((z+1)/(z-1))$ along $[-1,1]$, and the plot on the right
shows that the convergence
as a function of $n$ closely matches that of Gauss quadrature.
Here as in our later plots, convergence eventually stagnates for reasons
related to rounding errors, which we do not attempt to analyze.

For a variant problem that Carl Gauss
could not have handled, we start from the observation
that although Gauss quadrature is commonly regarded as optimal,
this is only true in specialized senses.
Gauss quadrature is the unique $n$-point quadrature formula
that integrates polynomials of degree $\le 2n-1$ exactly, and relatedly, it is
optimal for integrands analytic in a Bernstein ellipse, which is a natural
domain for polynomials.
For integration of functions known to be analytic in
domains other than ellipses, however, it is suboptimal \cite{bakh,exactness}.
A more natural assumption to make about an integrand may sometimes
be that it is analytic
in an $\varepsilon$-neighborhood of $[-1,1]$, and in~\cite{ht}, quadrature
formulas were constructed for such classes based on conformal transplantation
of a Bernstein ellipse.  With rational approximation, we can get the same effect
without the conformal mapping
by changing the definition of the approximation set $Z$
above from an ellipse to a stadium, like this:
\medskip

{\small
\begin{verbatim}
        ep = 1/sqrt(20);
        Z = [-1i*ep+linspace(-1,1,100)'; 1+ep*exp(1i*pi*(-49:49)'/100)];
        Z = [Z; -Z];
\end{verbatim}
\par}
\medskip

\noindent Figure~\ref{gaussvfig} shows that for the same integrand
$f(z) = 1/(1+20z^2)$ as before, the convergence is now faster by a
factor of about $\pi/2$.
Compare Figure 3.4 of \cite{ht} and Figure 4.1 of \cite{exactness}.  Perhaps
more importantly, the quadrature nodes are more evenly spaced, a property that
can be advantageous in time-stepping applications when there are
CFL stability constraints~\cite{kte}.
The quadrature weights, which remain positive apart from complex
perturbations, are also closer to uniform.\footnote{The reader may
be puzzled how one can have nearly equispaced nodes without coming up
against the Runge phenomenon of exponentially large weights with alternating signs.
The explanation is that that is a phenomenon of polynomial interpolation, whereas
here, the quadrature rule is implicitly interpolating by rational functions \cite{ht}.}

\section{Application 2: Quadrature with nearby singularities}

\begin{figure}
\begin{center}
\includegraphics[trim=0 125 20 7, clip, scale=.88]{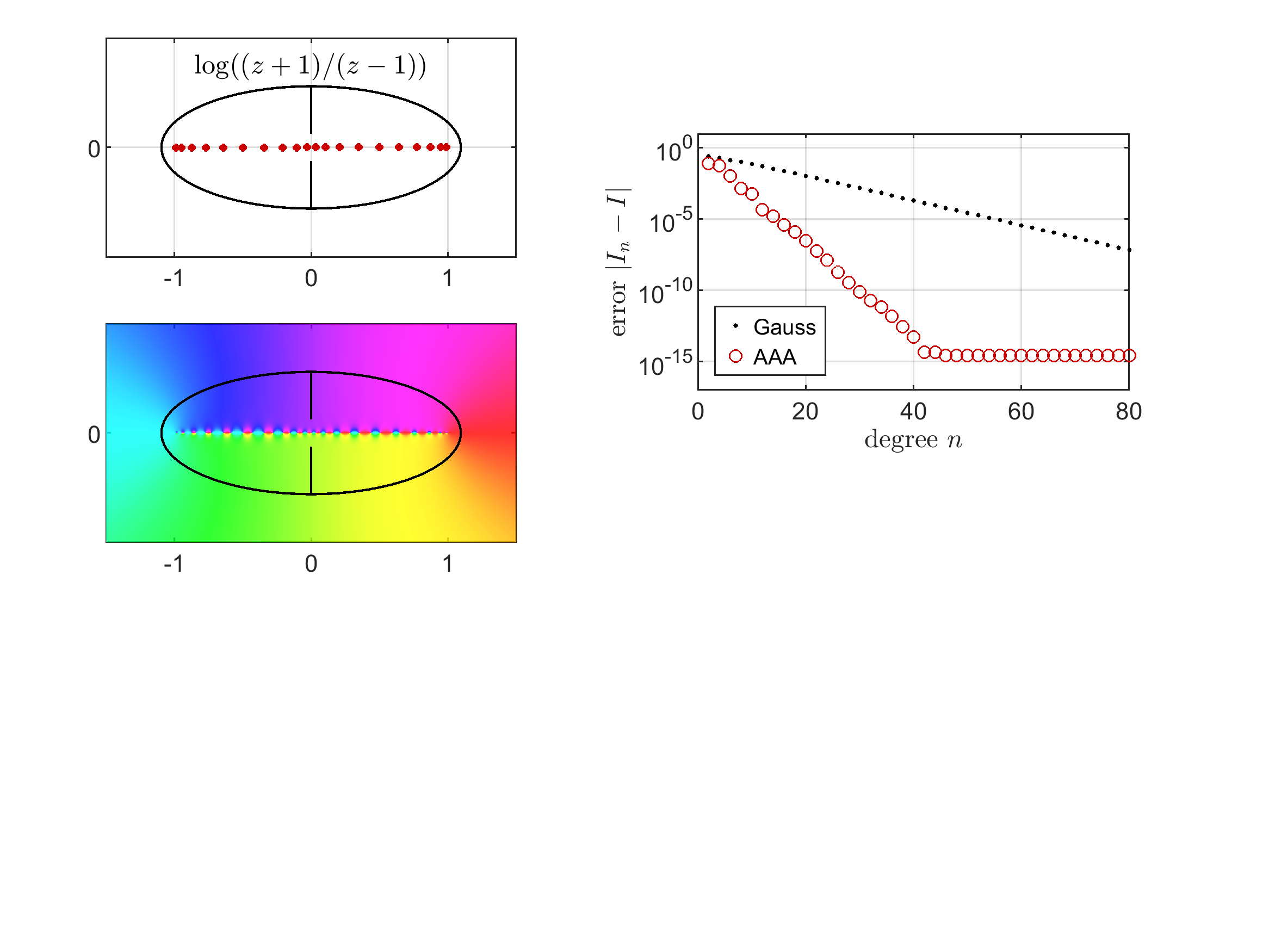}
\vspace*{-5pt}
\caption{\label{singfig}Construction of a quadrature
formula for an integrand with a conjugate pair of
singularities near $[-1,1]$, motivated by {\rm\cite{tee}}.  The test integrand
is $1/(1+100 z^2)$.}
\end{center}
\end{figure}

A theme of this article is that transplantation by conformal maps
has been a tool often used for enhancing accuracy of integrals.
A remarkable success in this direction was due to Wynn Tee in his
construction of numerical methods for time-dependent blow-up and
shock PDE problems \cite{tee}.  (The second author was
the junior partner in this work.)  Using an adapted spectral spatial
grid with as few as 56 points, Tee was able to resolve the solution
to $u_t^{} = u_{xx}^{} + e^u$ on $[-1,1]$ with homogeneous initial
and boundary conditions up to a time within $10^{-8}$ of the blowup
time of $t_c^{}\approx 3.54466459$.  Tee's method can be described
with reference to Figure~\ref{singfig}.  The difficulty in these
problems is a pair of complex singularities approaching the real axis
as $t\to t_c^{}$, for singularities at distance $\varepsilon$ would
ordinarily require a grid of $O(\varepsilon^{-1})$ points.  If the
singularities are regarded as the ends of branch cut slits in the
complex plane, however, he showed that a conformal transplantation
can be used to send them far away, improving $O(\varepsilon^{-1})$
to $O(|\log \varepsilon|)$.

All this required careful analysis of conformal maps, but in the
experiment shown in the figure, we have achieved the same effect (for
quadrature rather than PDE\kern .3pt s) by rational approximation.
The same Cauchy transform $\C(s) = \log((z+1)/(z-1))$ as in the last section is
now approximated on a Bernstein ellipse to which a
conjugate pair of slits extending to $\pm 0.1 i$ have been added.  For
the test integrand $f(z) = 1/(1+100 z^2)$, whose singularities
lie at the ends of the slits, one sees about a five-fold speedup over
Gauss--Legendre quadrature.  Closer singularities would amplify
the effect.  In Tee's paper he estimates the singularity positions
on the fly with Chebyshev--Pad\'e approximation, which we would now
propose to replace by AAA.

\begin{figure}
\begin{center}
\includegraphics[trim=0 125 20 3, clip, scale=.88]{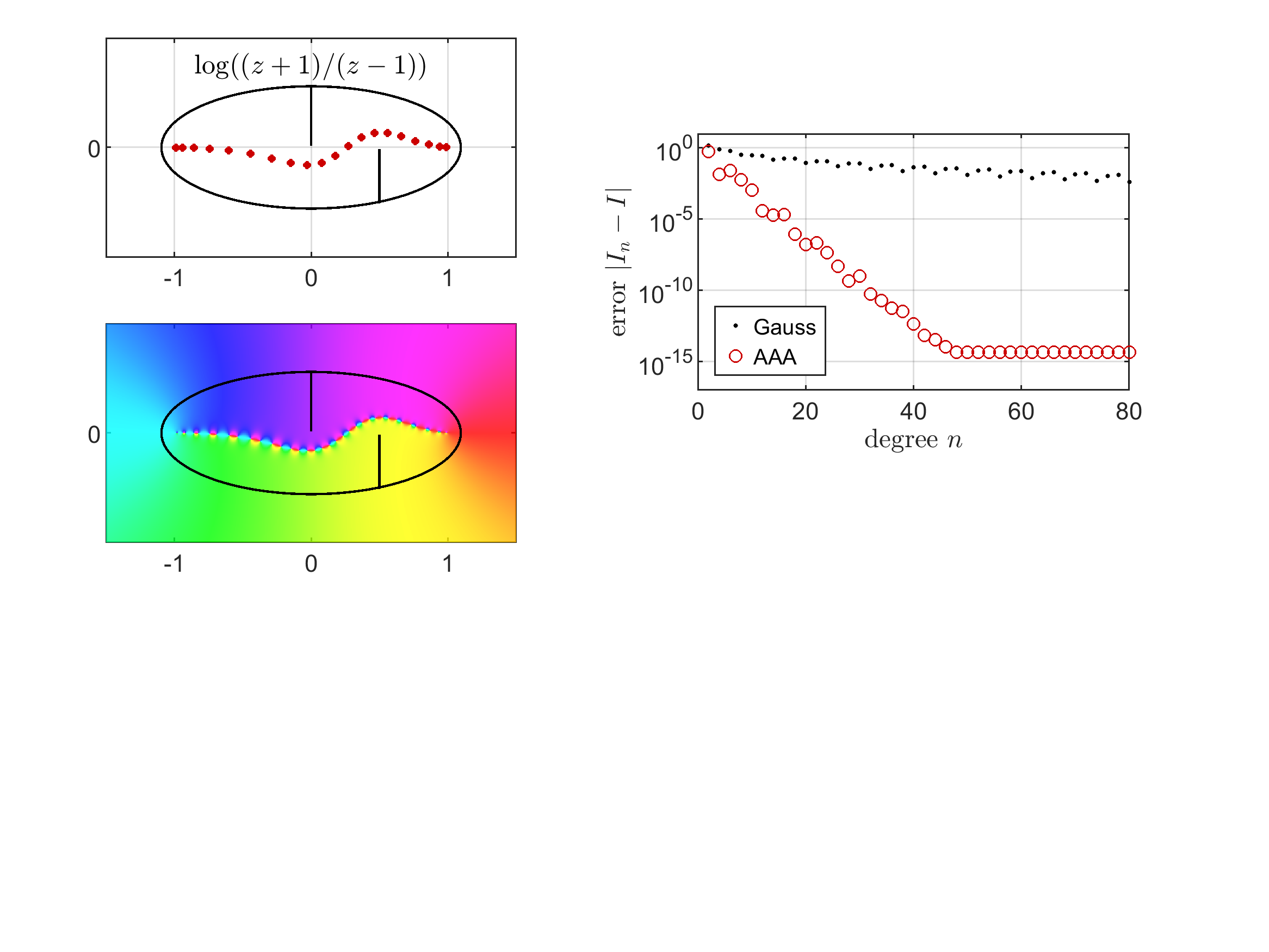}
\vspace*{-5pt}
\caption{\label{singvfig}Variant of Figure {\rm
\ref{singfig}} inspired by {\rm\cite{haletee}}
in which there may be several singularities
near $[-1,1]$, and not just conjugate pairs.  The
test integrand is
$\log([z-.01i]i) + \log([z-(.5-.01i)]/i)$.}
\end{center}
\end{figure}

After the completion of \cite{tee}, Hale and Tee generalized
Tee's method to problems with multiple conjugate pairs of
singularities near the approximation interval \cite{haletee}.
Their generalization was again highly successful, but required
expert Schwarz--Christoffel mapping.  As a variant problem for this
section, Figure \ref{singvfig} shows the same kind of improvement as in
Figure~\ref{singfig} for
a problem with several singularities.  To highlight the flexibility
of the method we have chosen a complex integrand whose singularities
do not lie in conjugate pairs, namely $f(z) = \log([z-.01i]i) +
\log([z-(.5-.01i)]/i)$.  The snaky curve in the figure emphasizes
the fact that in general, rational approximation chooses the contour
as well as the nodes and weights.  The examples of
Figures \ref{fig:osc}, \ref{fig:oscv},
\ref{talbotfig}, \ref{talbotvfig}, \ref{circlevfig},
\ref{funAfig}, and \ref{funAvfig}
are also notably of this kind.  Indeed, near-optimal
rational approximations are always
choosing a near-optimal contour, and this merely fails to be
apparent when the contour is simply a real interval for reasons
of symmetry.

\section{Application 3: Singular integrals}

\noindent
The method of section~2 is straightforward to adapt to a variety of challenging
(and interesting!) non-constant weight functions $w(z)$ in
(\ref{defn}). The new ingredient is that we
need samples from the Cauchy transform of $w(z)$ to feed into AAA.\ \ In
some special cases, just as for $w(z)=1$ in the
last two sections, the Cauchy transform can
be determined analytically. In practice, it is convenient to
evaluate it numerically on the contour
$\Gamma$ using an adaptive quadrature rule. In this and the next
section we
demonstrate this approach in two settings: first weakly
singular integrals, then highly oscillatory ones.

For the first example, suppose $w(z)$ in (\ref{defn})
is integrable but not necessarily smooth or bounded.
Numerical
methods are often based on adaptive node refinement, analytic
expansions, domain mapping, or hybrid schemes that combine these
elements (for example, see ``quadrature by expansion''~\cite{klockner} and
references therein).  Adaptive refinement may evaluate $f(z)$ many
times, while analytic expansions and domain mapping techniques, for
maximal efficacy, are derived on a case-by-case basis for different
classes of weight functions.  However, rational approximation
offers the alternative of easy construction of Gauss-like formulas for
arbitrary weights.
To illustrate, consider the asymmetric Jacobi weight
$w(z) = (1+z)^{3/2} (1-z)^{-1/2}$.  This weight has
algebraic singularities at both endpoints and blows up at the right
endpoint.
To obtain a quadrature rule with AAA, we compute samples of
$\mathcal{C}(s)$ along the Bernstein ellipse of section~3 with 
Gauss--Kronrod quadrature as implemented in {\tt quadgk} in MATLAB:
\medskip

{\small
\begin{verbatim}
    rho = 1/sqrt(20) + sqrt(21/20);
    c = rho*exp(2i*pi*(1:400)/400);
    Z = (c+1./c)/2;
    w = sqrt(1+z).^3./sqrt(1-z);
    F = @(z) arrayfun(@(zz) quadgk(@(z) w(z)./(zz-z),-1,1), z);
    [r,pol,res] = aaa(F,Z,'degree',20,'sign',1);
\end{verbatim}
\par}
\medskip

\begin{figure}
\begin{center}
\includegraphics[trim=0 4.8cm 0 0,clip,scale=.82]{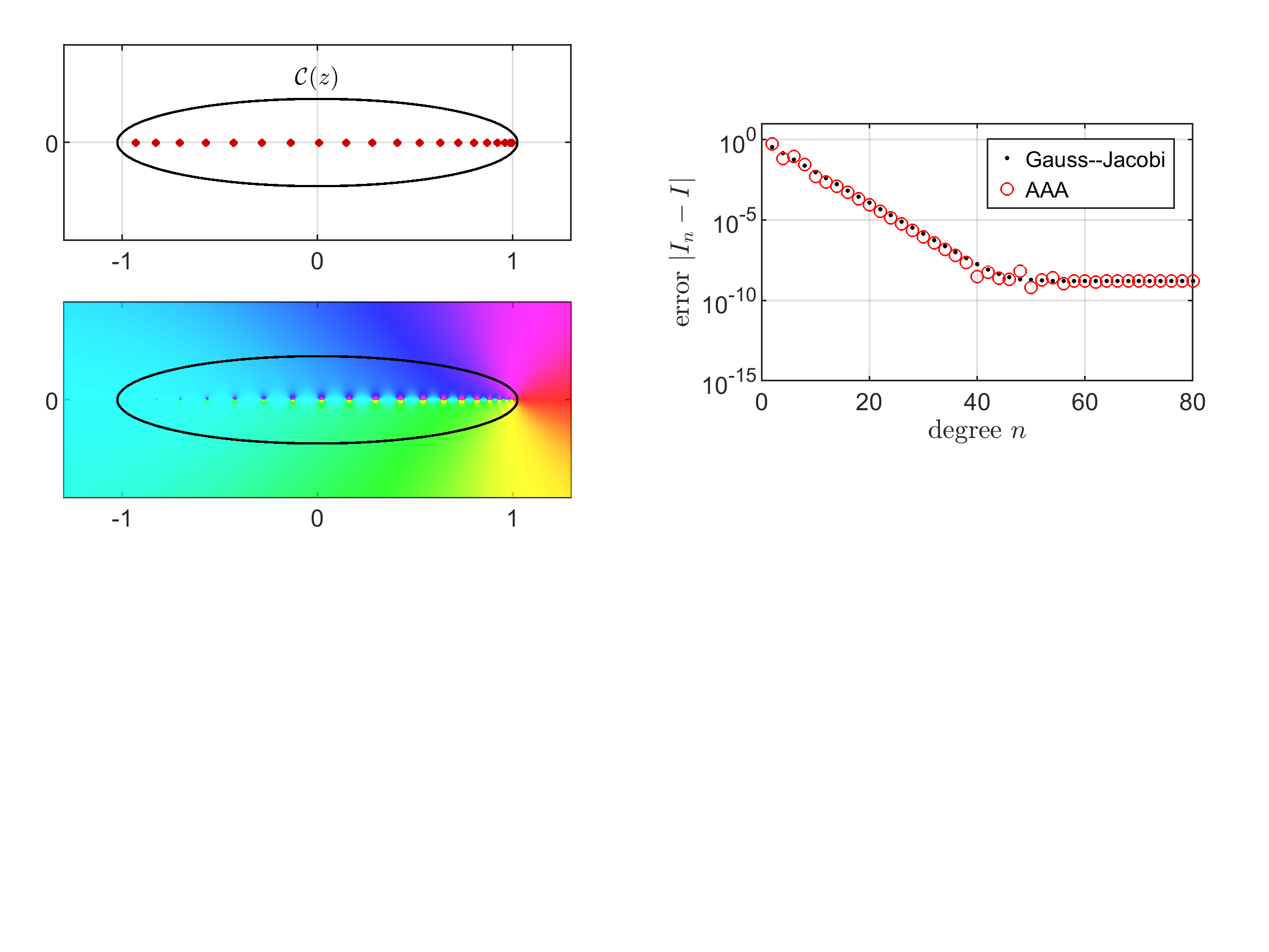}
\end{center}
\vspace*{-4pt}
\caption{\label{fig:sing} AAA quadrature for the singular Jacobi weight
function $w(z) = (1+z)^{3/2}(1-z)^{-1/2}$.  The AAA poles cluster asymmetrically
to capture the blow-up at $x=1$.}
\end{figure}

\noindent
For higher-degree
rules, one should adjust the default tolerances
in {\tt quadgk} to sample $\C(s)$ with sufficient
accuracy. For this section, we set both \verb|'AbsTol'| and
\verb|'RelTol'| to $10^{-13}$. The left panel of Figure~\ref{fig:sing}
displays the quadrature nodes (poles of 
$\rn(s)\approx \C(s)$) along with the phase portrait of $\rn$.
The quadrature nodes
cluster asymmetrically toward the right endpoint to capture
the blow-up of $w(z)$ there,
and the phase portrait makes it clear that that is where most of the action is.
The Gauss--Jacobi nodes and weights
associated with this weight function exhibit similar behavior (not shown).
The right panel compares the convergence of AAA
quadrature, for the same function $f(z)$ as in section~3, with that
of the Gauss--Jacobi quadrature rule computed with \texttt{jacpts}
in Chebfun.
Clearly rational approximation has discovered rules that are essentially
as effective as Gauss quadrature.

\begin{figure}
\begin{center}
\includegraphics[trim=0 4.8cm 0 0,clip,scale=.82]{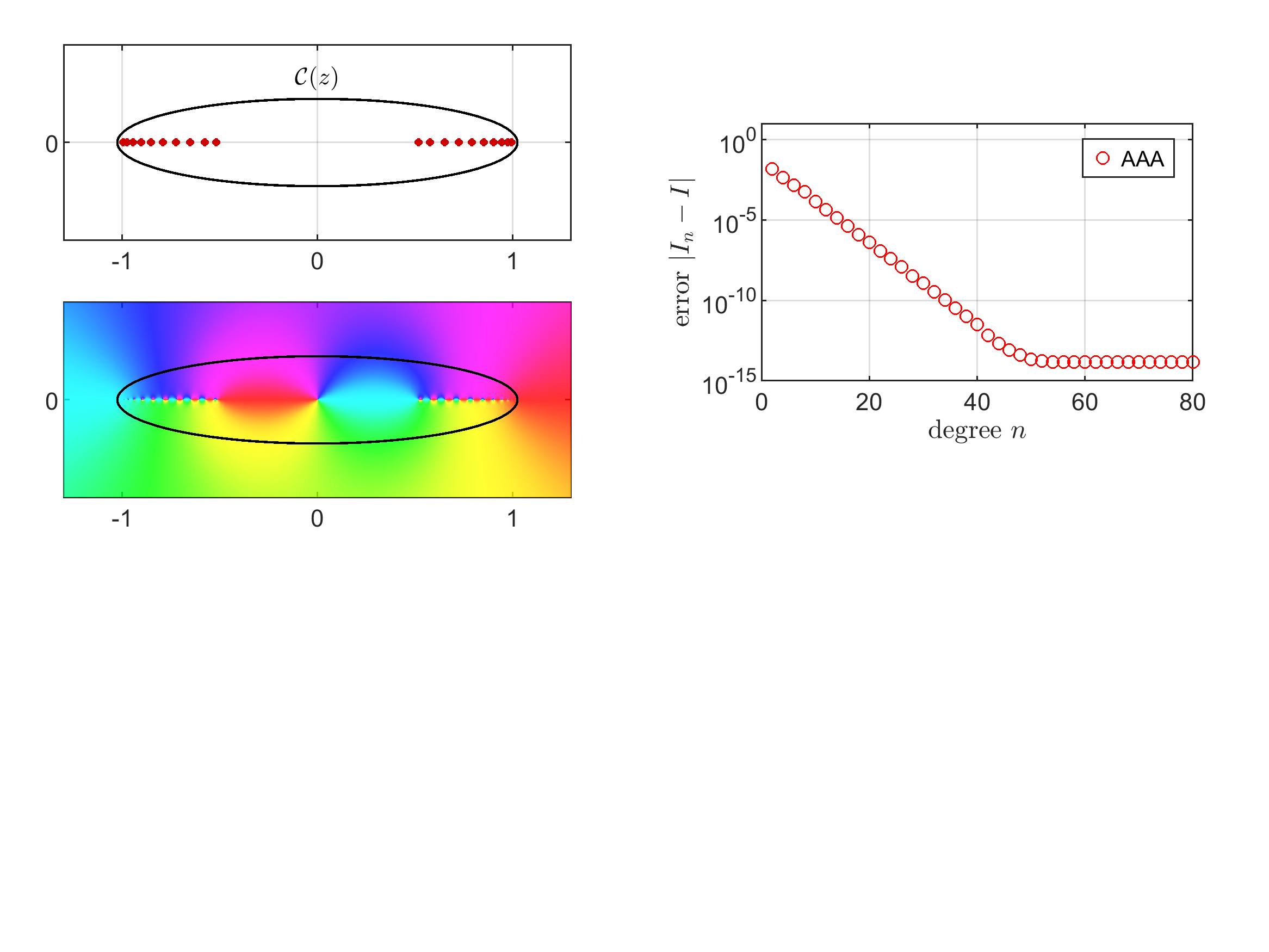}
\end{center}
\vspace*{-4pt}
\caption{\label{fig:singv} AAA quadrature for the
weight function {\rm (\ref{oscwt2})}.  Poles appear only in the portions of $[-1,1]$
where $w(z)$ is nonzero.}
\end{figure}

Our variant for this section concerns the non-classical weight function
on $[-1,1]$ defined by
\begin{equation}
w(z) = \begin{cases} \sqrt{1-z^2} & 0.5\le |z|\le 1, \\
0 & \hbox{otherwise.} \end{cases}
\label{oscwt2}
\end{equation}
This weight function is supported on the disjoint intervals
$[-1,-0.5\kern .4pt ]$ and $[\kern .4pt 0.5,1]$ and is not differentiable
at any of the endpoints.  Weight functions supported on
disjoint intervals
present a common challenge in the computation of
response integrals in quantum chemistry and 
materials science, where points of non-differentiability at
``band edges'' or in the ``bulk bands''
are often associated with fascinating physical phenomena.  
High-order quadrature schemes must take care to resolve such singularities, which can 
degrade the convergence of popular ``out-of-the-box'' methods like the Kernel
Polynomial Method (KPM)~\cite{silver,weisse}.

AAA generates a new quadrature rule for the non-classical weight
function with a single modification to the code above,
namely, defining the new weight function.  Figure~\ref{fig:singv} is the
analogue of Figure~\ref{fig:sing} for the new AAA quadrature rule, and
one can see that the quadrature nodes are confined to the support of
the weight function. For
this example, one could also design a specialized rule
by partitioning the integration domain into
$[-1,-0.5]$ and $[\kern .4pt 0.5,1]$ and applying different Gauss--Jacobi
rules on each subdomain. However, this strategy requires detailed
knowledge of the location and type of each singularity. The AAA
approach achieves similar performance while only requiring evaluation
of the weight function on $[-1,1]$.

\section{Application 4: Oscillatory integrals}

For our next application,
suppose we want to compute an \textit{oscillatory} integral of the form 
\begin{equation}\label{eqn:osc}
I = \int_{-1}^1 f(z) \exp(i \kern .4pt \omega g(z)) \kern 1pt dz,
\end{equation}
where $g(z)$ is a smooth {\em phase function} and $\omega$ is a
real frequency. Such oscillatory integrals arise
in scattering problems and related applications in, e.g.,
imaging. Traditional Gauss quadratures require many samples of
$f(z)$ when $\omega$ is large, as they must satisfy the Nyquist
sampling criterion associated with polynomial interpolation of
the integrand. More efficient methods, which improve as
$\omega$ increases, rely on specialized asymptotic
expansions that exploit oscillation and cancellation in the
integrand~\cite{dhibook,filon,levin82}.
When $g(z)$ is analytic, the key ideas are those
of stationary phase, to identify points of dominant contribution,
and steepest descent, to identify paths in the complex plane
along which the integrand decays exponentially.  By deforming
the contour off of $[-1,1]$ along paths of steepest descent
through points of stationary phase, one can truncate and
apply Gauss quadrature on the deformed integral
to derive methods that are convergent and of optimal order as
$\omega\rightarrow\infty$~\cite{deano09,huybrechs06}.

Remarkably, the AAA approximation to the Cauchy transform of the
weight function $w(z) = \exp(i\kern .4pt\omega g(z))$ appears
to generate quadrature rules of similar structure
and similar efficiency.
To illustrate, fix $\omega
= 25\pi \approx 78.53$ and consider the simplest phase function,
$g(z) = z$, corresponding to a band-limited Fourier transform. Since
$g'(z) = 1$ is non-vanishing, there are no stationary points, and the
dominant contribution to $I$ comes from the endpoints. The classical
``two-point'' quadratures associated with highly oscillatory integrals
approximate $I$ with a weighted combination of $f$ and its
derivatives evaluated at the endpoints.  For a fresh account with
modern developments, see~\cite{iserles06}.

\begin{figure}
\begin{center}
\includegraphics[trim=30 4cm 30 0,clip,scale=.89]{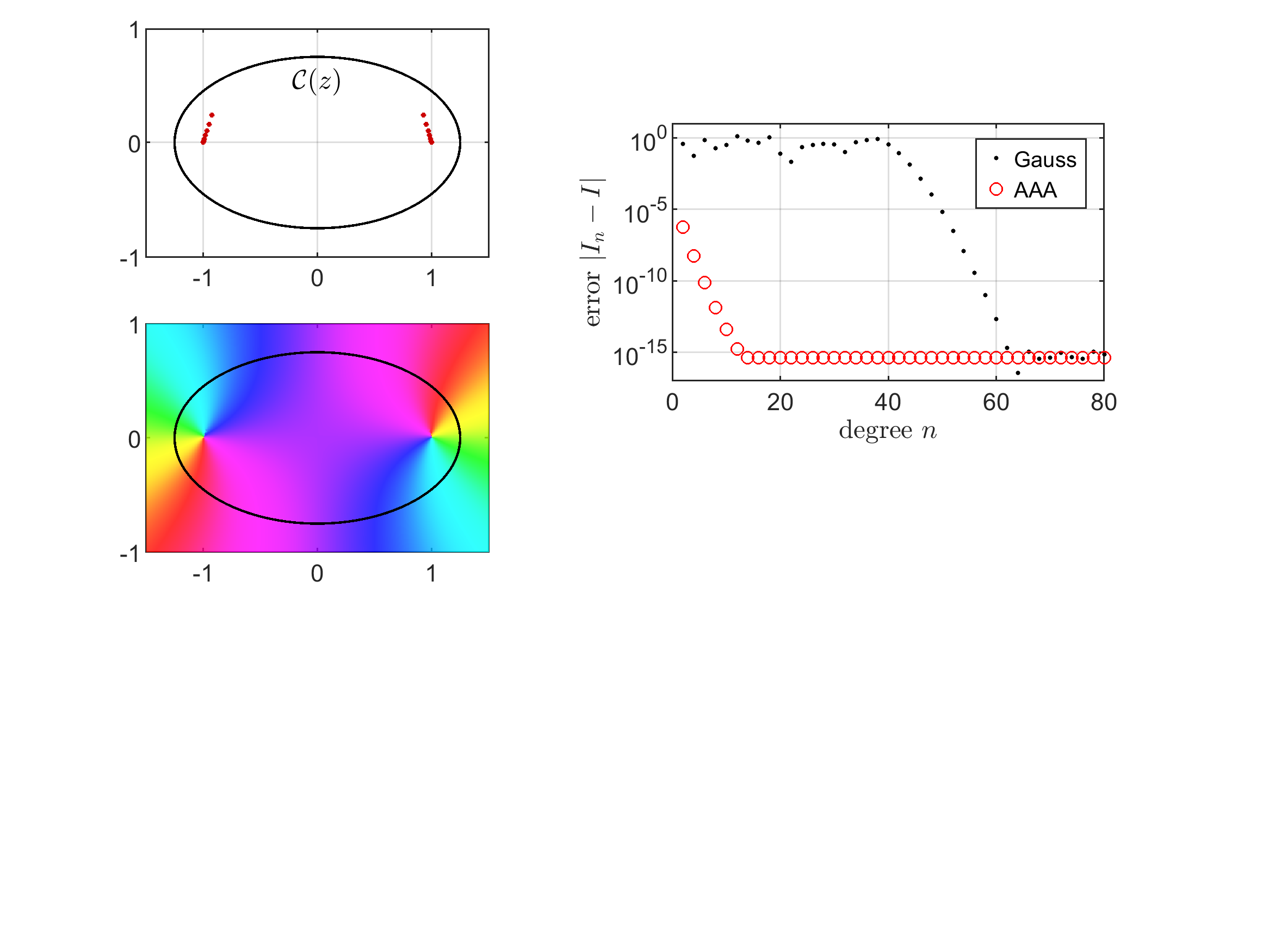}
\end{center}
\vspace*{-9pt}
\caption{\label{fig:osc} AAA quadrature for the oscillatory integral
in $(\ref{eqn:osc})$ with $\omega=25\pi$ and $g(z) = z$. The quadrature
nodes cluster near the endpoints along arcs that are close to
paths of steepest descent. The right panel shows convergence
for $f(z) = 1/(1+z^2/4)$ as the degree of
the quadrature rule is refined and compares with Gauss--Legendre quadrature.}
\end{figure}

\begin{figure}
\begin{center}
\includegraphics[trim=30 4cm 30 0,clip,scale=.89]{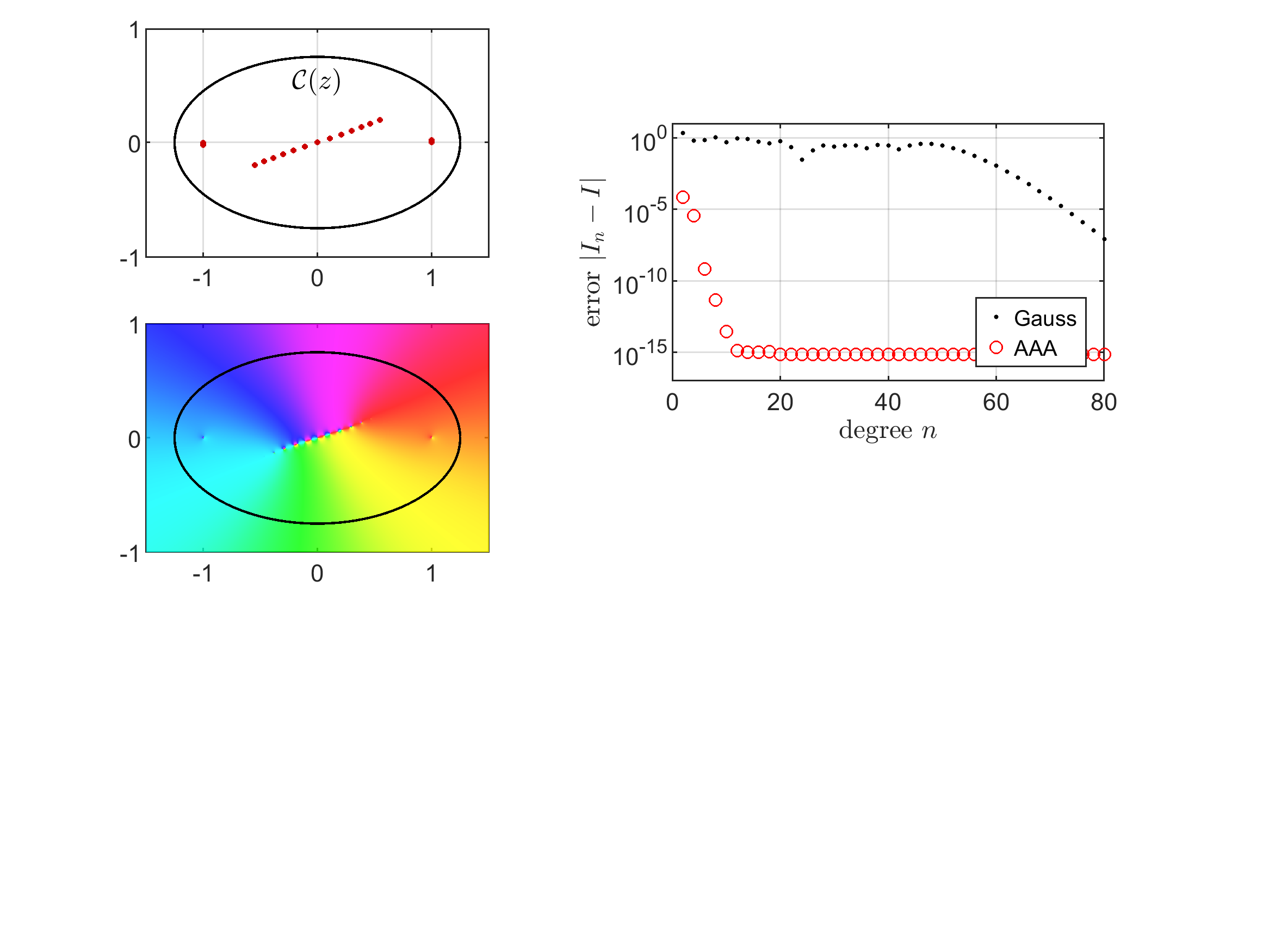}
\end{center}
\vspace*{-9pt}
\caption{\label{fig:oscv} Same as in Figure~$\ref{fig:osc}$ but for
the weight function $g(z) = z^4$.}
\end{figure}

Now, as in the previous section, we compute the Cauchy
transform of $w(z) = \exp(25\pi i z)$ on a Bernstein ellipse
with $\rho=2$ using {\tt quadgk}. The samples of the Cauchy
transform are fed to AAA, which generates a rational approximation of
user-specified degree, whose poles and residues become the quadrature
nodes and weights. Figure~\ref{fig:osc} shows that these quadrature
nodes cluster densely at the endpoints of the interval. Moreover,
they lie along nearly vertical arcs in the upper half-plane, with a
slight inclination away from the ellipse. Further
experiments (not shown) reveal that when the Cauchy transform is
sampled on a smaller Bernstein ellipse, the nodes crowd inward
as if confined to the interior of the ellipse. When the ellipse
is made larger, the nodes tend toward vertical lines in the upper
half-plane---the paths of steepest descent for this Fourier integral.


For the variant, we consider the phase function $g(z) = z^4$,
which has a stationary point of order $3$ at $z=0$. The associated
path of steepest descent slices through the origin in the complex
plane at an angle of $\pi/8$. The quadrature nodes computed by AAA
again cluster at the endpoints, but now some of them distribute
themselves evenly along an arc close to the path of steepest descent
through $z=0$. When the Cauchy transform is sampled on a smaller
Bernstein ellipse (not shown), the nodes distribute themselves
along an arc that is close to the path of steepest descent near
the real line but curves away to fit within the ellipse off
of the real line.

\section{Application 5: Hankel contours for inverse Laplace transforms}

John Butcher in the 1950\kern .6pt s and Alan Talbot in the
1970\kern .6pt s introduced the powerful technique of computing
inverse Laplace transforms via discretizations of
integrals over so-called {\em Hankel contours,} 
\begin{equation}
I = \int_\gamma \kern -2pt e^z f(z) \kern 1pt dz 
\,\approx\, \sum_{k=1}^n c_k f(z_k).
\label{hankel}
\end{equation}
A Hankel contour $\gamma$ is an infinite arc in the complex
plane wrapping around $\negreal$ from $-\infty-0\kern .4pt i$ to
$-\infty+0\kern .4pt i$ (and thus not precisely a homeomorphism
of $[-1,1]$).  Here $f$ is a function analytic in the simplest
cases in $\complex\backslash \negreal$ which may be a scalar for
applications to special functions or a large-dimensional matrix
function for applications in computational science.  Many researchers
have worked with such techniques including Gallopoulos, Gavrilyuk,
Gil, L\'opez-Fern\'andez, Lubich, Luke, Makarov, McLean, Palencia,
Sch\"adle, Sheen, Sloan, Saad, Schmelzer, Temme, Thom\'ee, and
Trefethen, and the one who has investigated them most systematically
is Weideman \cite{talbot,bromwich}.  See these papers and \cite{trap}
for technical details and references.

To address (\ref{hankel}) by rational approximation, we interpret
$e^z$ as a weight function on some Hankel contour; exactly which
one doesn't matter.  Since $e^z$ is analytic, it follows from
(\ref{cauchytrans}) that the Cauchy transform $\C(s)$ is $e^s$
itself.  The contour $\gamma$ is deformed inward to
$\Gamma =\negreal$ (traversed from left to right along the bottom and
then right to left along the top), and then in principle~$\Gamma$
should also include an outer boundary on which the value to be
matched will be~$0$.  However, this part of $\Gamma$ does not need
to be considered explicitly, as a rational approximation to $e^z$
on $\negreal$ will necessarily be $\approx 0$ for $z\to\infty$.
We can compute one with the code

\medskip

{\small
\begin{verbatim}
        Z = -logspace(-3,4,300)';
        F = exp(Z);
        [r,pol,res] = aaa(F,Z,'degree',14); 
\end{verbatim}
\par}
\medskip

\begin{figure}
\begin{center}
\includegraphics[trim=30 90 20 3, clip, scale=.88]{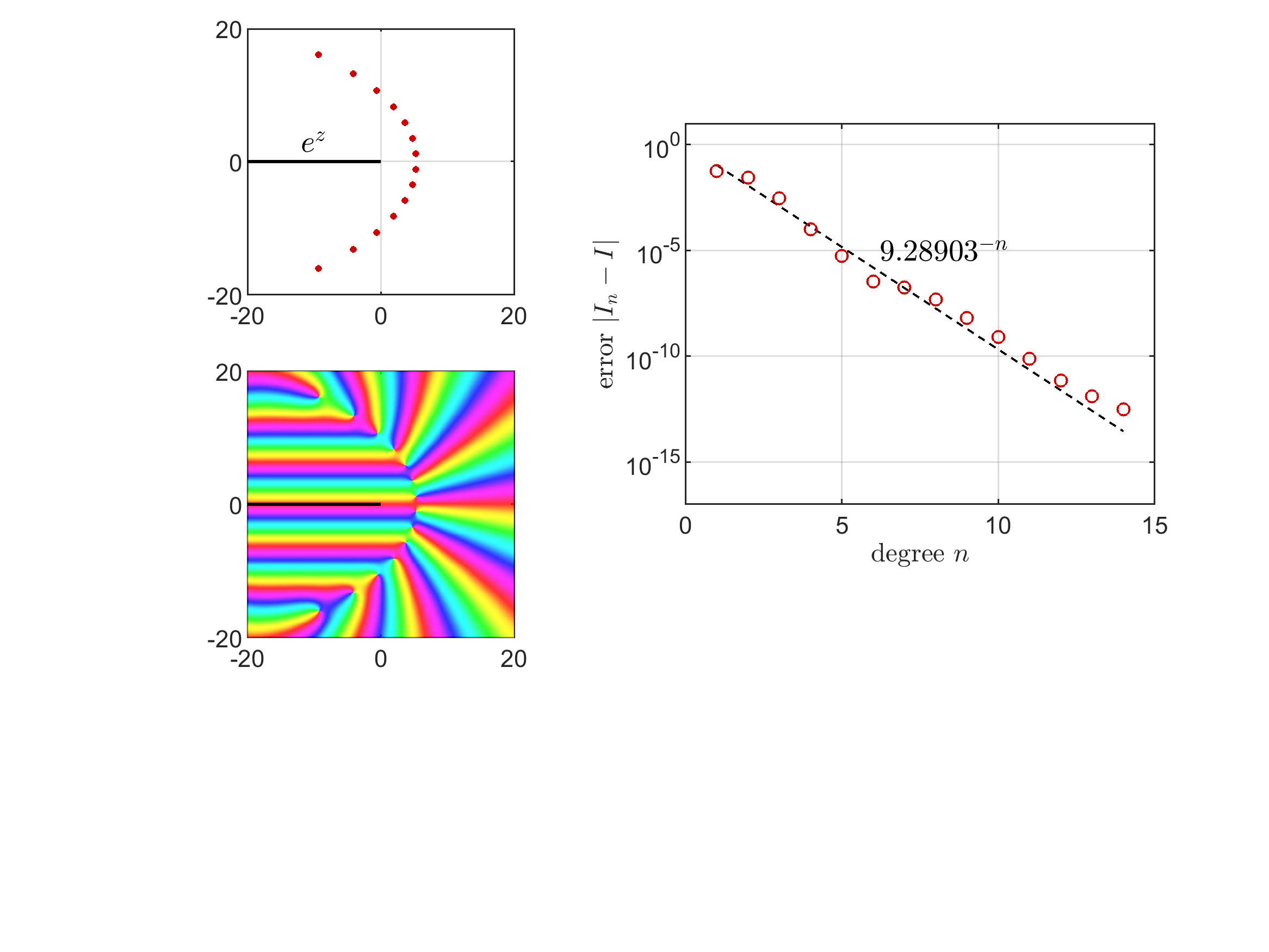}
\end{center}
\caption{\label{talbotfig}
Evaluation of the inverse Laplace transform by rational approximation following
Butcher, Talbot, and many later authors.  On
the left, quadrature poles with $n=14$, giving $12$-digit accuracy for a test
problem based on just $7$ matrix solves (if real symmetry is exploited).
The stripes in the phase portrait match those of\/ $e^z$ to the left of the curve of poles.
On the right, convergence as a function of $n$, matching the
optimal rate for best approximation.}
\end{figure}

\noindent The degree $14$ is fixed as the highest possible
before reaching machine precision \cite[chap.~20]{atap}.
Figure~\ref{talbotfig} shows the very satisfactory result, with poles
closely matching those in the papers cited above.  A numerical calculation based on
another integrand with integral 1,
\medskip

{\small
\begin{verbatim}
        f = @(z) -exp(1)./(1+z);
        In = res.'*f(pol)
\end{verbatim}
\par}
\medskip

\noindent gives the result
$0.99999999999937 + 0.00000000000003\kern .4pt i$
with error
$6.3\cdot 10^{-13}$.  The convergence as a function of $n\to\infty$
closely tracks the optimal rate $(9.28903\dots)^{-n}$ \cite[Thm.~25.2]{atap}.

The elegant curve of poles on the left in Figure~\ref{talbotfig}
echoes similar images to be found in Figures 3.1--3.3 and 4.3
of \cite{talbot}, where quadrature formulas are compared based on
parabolic, hyperbolic, and cotangent (``Talbot'') contours as well as
true best approximations of $e^z$ on $\negreal$.  In the first three
cases, exponential convergence rates on the orders of $2.85^{-n}$,
$3.20^{-n}$, and $3.89^{-n}$ are obtained, whereas as we see in
Figure~\ref{talbotfig}, rational approximation gives the optimal rate
of about $9.29^{-n}$.  This approximate doubling of the convergence
rate, which we think of as the ``Gauss quadrature factor of 2,''
arises in many contexts of rational approximation in the difference
between estimates obtained from potential theory arguments assuming
$n+1$ interpolation points of a rational approximation and the
improved estimates, related to orthogonality, that come when there
are $2n+1$ interpolation points.  See the end of \cite{notices}
for a brief discussion and \cite{rakh} for theorems.

\begin{figure}
\begin{center}
\includegraphics[trim=30 200 20 3, clip, scale=.92]{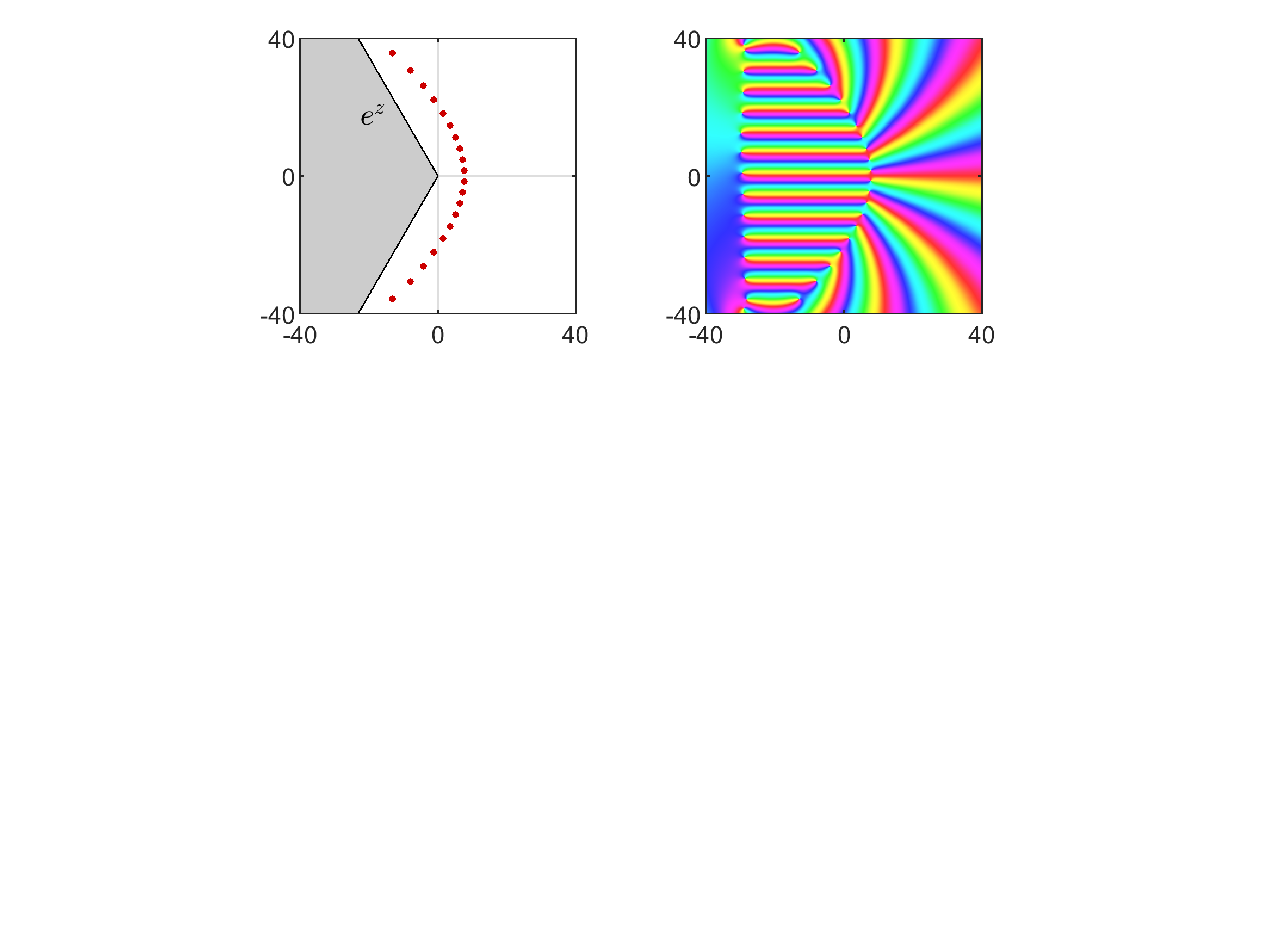}
\end{center}
\caption{\label{talbotvfig}Variant problem for a sectorial operator
as in\/ {\rm \cite{lps}}.
The reason the stripes terminate at around real part $-35$ is
that $e^{-35}$ is on the order of the accuracy of this approximant, close
to machine precision.}
\end{figure}

Our variant problem for this section concerns a setting in
which it is supposed that the function $f$, which may be derived
from a matrix
or operator, is {\em sectorial,} meaning that it is analytic
and satisfies a growth condition outside a wedge in the left half-plane \cite{lps,bromwich}.
A small change in the definition of the approximation domain,
\medskip

{\small
\begin{verbatim}
        Z = -logspace(-3,4,300)';
        Z = [flipud(Z)/exp(.333i*pi); Z*exp(.333i*pi)];
\end{verbatim}
\par}
\medskip

\noindent gives the images of
Figure~\ref{talbotvfig}, now with the degree increased from $n=14$
to $n=20$ since the problem is slightly harder.

\section{AAA approximation of sign functions}

The last five sections could have been written five years ago,
based on the AAA and AAA-Lawson algorithms of \cite{aaa} and
\cite{lawson}.  In the next two sections we will turn to examples
involving approximation on two disjoint contours (see the right
image of Figure~\ref{schematic}), where a rational
approximation is needed to a function with two branches, and here,
these original algorithms proved unreliable.  However, modifications
were introduced in 2024 that make problems of this kind tractable
too~\cite{zolo}.  These variants are invoked by specifying the
\verb|'sign'| and/or \verb|'damping'| options in calls to {\tt aaa}.
The first of these adjustments bases the AAA choice of a barycentric
weight vector on a weighted blend of all the singular vectors of
a Loewner matrix rather than just one, and the second introduces a
damping factor in the update formula of the AAA-Lawson iteration.
We will not discuss the details here, as they are treated
in~\cite{zolo} and the topic of
this paper is the link from rational approximation to quadrature,
not the mechanics of rational approximation.

\begin{figure}
\begin{center}
\includegraphics[trim=55 30 55 10, clip,scale=0.5]{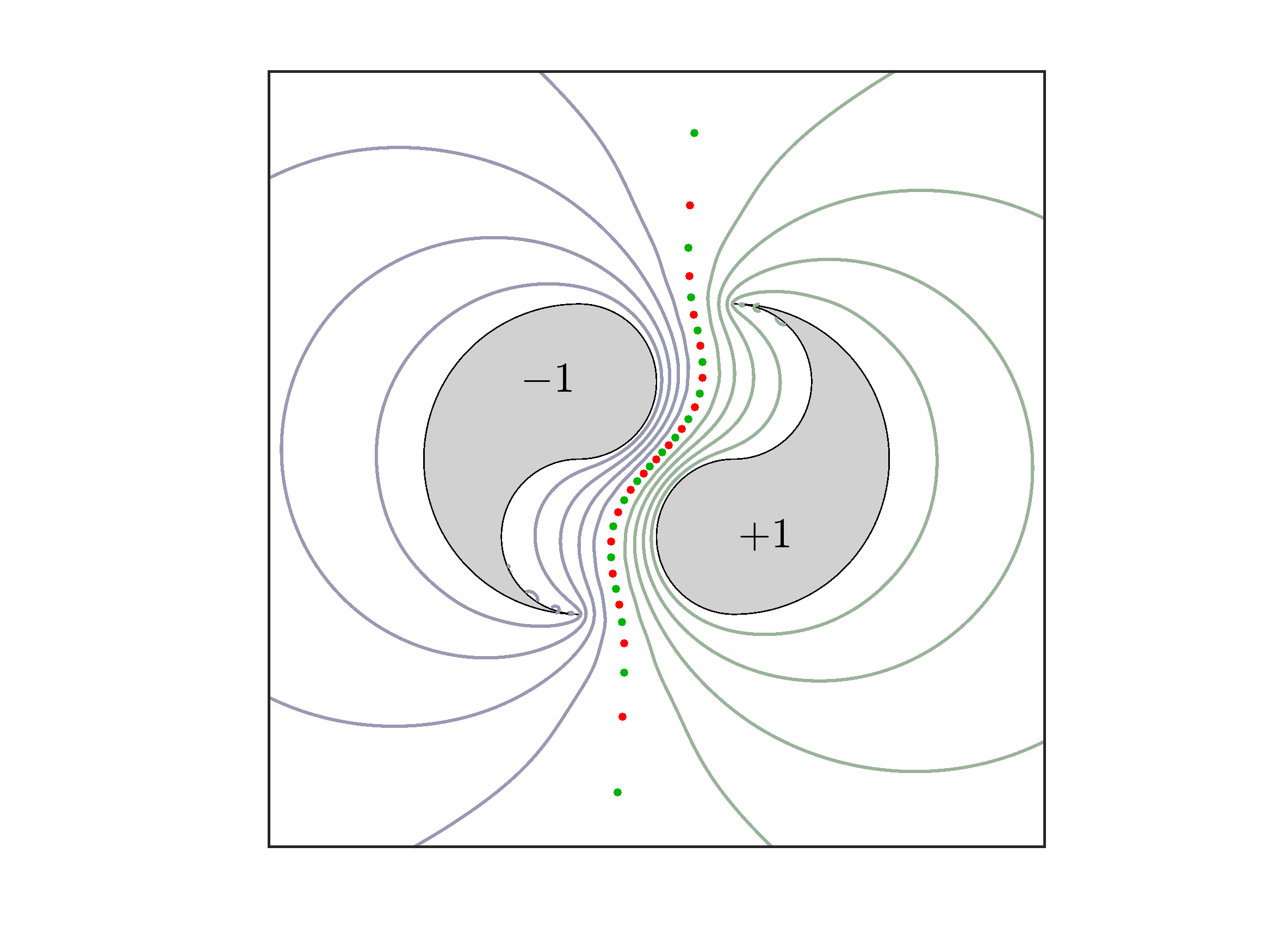}
\includegraphics[trim=55 30 55 20, clip,scale=0.5]{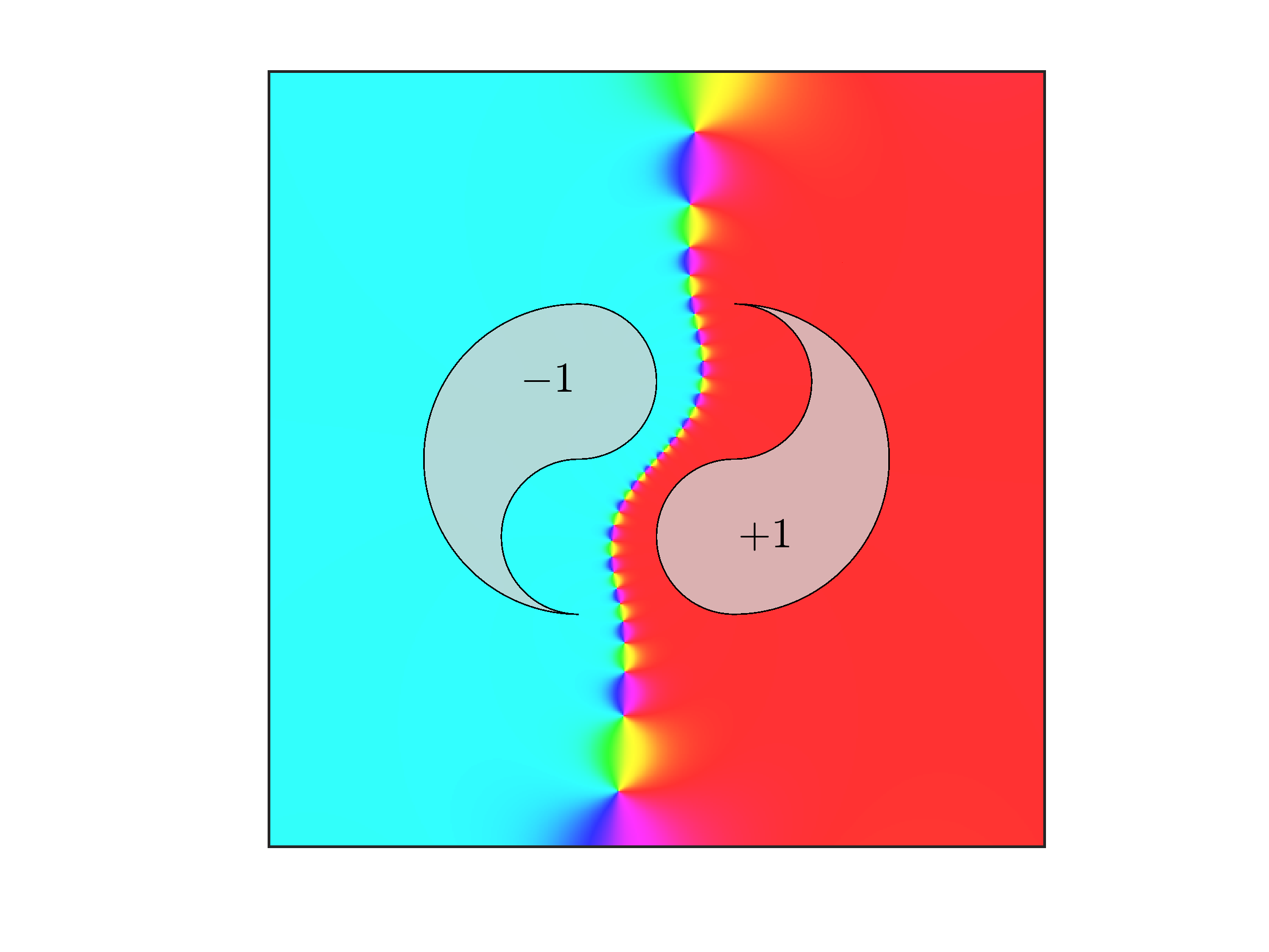}
\end{center}
\vspace*{-3pt}
\caption{\label{yinyangfig} Degree\/ $20$ rational approximation $r(z)$ of a sign function on a 
two-component yin-yang domain in the complex plane.
On the left, the domain is shown
together with poles (red) and zeros (green) of the rational approximation,
defining an approximate branch cut, as well as error contours 
$|r(z) - (\pm 1)| = 10^{-1},\dots,10^{-4}$.  On the right, a phase portrait
of\/ $r$.}
\end{figure}

To illustrate the new capabilities,
here is an example adapted from~\cite{zolo} of AAA computation
of an approximate sign function.  Suppose it is desired to find
a rational function $\rn$ of degree $n=20$ that comes as close as
possible to taking the values $-1$ and $+1$ on the two components of
the yin-yang domain of Figure~\ref{yinyangfig}.
This is an example of ``Zolotarev's 4th problem,''
also known as the Zolotarev
sign problem \cite{it,zolo}.  We discretize the domain by 300 points
on the boundary of each component.  
\medskip

{\small
\begin{verbatim}
        c = exp(1i*pi*(1:100)'/100)/1i;
        yin = [-c; conj(1i*c)/2i-.5i; c/2+.5i] - 1/2;
        yang = -yin;
        Z = [yin; yang];
\end{verbatim}
\par}
\medskip

\noindent
No single analytic function can match the values
$-1$ and $1$ exactly; there must be a branch cut.  Near-best rational approximations
mimic this behavior by forming approximate branch cuts out of
strings of poles and zeros.
The rational approximation is computed in about $0.7$ s on a laptop
by these commands:
\medskip

{\small
\begin{verbatim}
        F = [-ones(size(yin)); ones(size(yang))];
        [r,pol,res,zer] = aaa(F,Z,'degree',20,'sign',1);
        plot(yin,'-k'), hold on, plot(yang,'-k')
        plot(zer,'.g'); plot(pol,'.r');
\end{verbatim}
\par}
\medskip

\noindent
The error contours in the figure can be understood at least
approximately---and exactly, in the limit $n\to\infty$---by methods
of potential theory~\cite{notices}.  The potential function in
question is generated by positive and negative charges on the two
components distributed in a minimal-energy configuration.

The example of Figure~\ref{yinyangfig} can be regarded as a prototype
of the kind of rational approximation needed for applications
involving Cauchy integrals over closed contours $\gamma$ that
separate $\complex$ into an interior and an exterior region.

\section{Application 6: Trapezoidal rule on a circle}

Here and in the next section we
derive quadrature formulas for problems with closed contours, beginning
with the simplest
case of an integral around a circle.
If $g$ is an analytic function on the unit circle
$\gamma = S = \{z\in\complex: |z|=1\}$, it has a Laurent series
\begin{equation}
g(z) = \sum_{k=-\infty}^\infty \ck z^k, \quad
\ck = {1\over 2\pi i} \oint_S  \! z^{-1-k} g(z) \,dz.
\end{equation}
Since the work of Lyness in the 1960\kern .4pt s, it has been
recognized that integrals like this 
can be accurately evaluated by
the trapezoidal rule with sample points at the $n$th roots of unity
$\{\zk\}$ for some $n$ \cite{lyness-moler,lyness-sande,trap}.
For any $f$ analytic on $S$,
\begin{equation}
{1\over 2\pi i} \oint_S f(z) \,dz
\approx \sum_{k=1}^n {\zk\over n} \, f(\zk) .
\label{circint}
\end{equation}
Indeed, this is the standard method for computing Taylor and Laurent
coefficients, with adjustments to optimize the radius of the circle
of evaluation depending on $f$~\cite{bornemann,fornberg}.  For a
review of the history of this technique with many more references,
see~\cite{akt}.  Applications in numerical linear algebra include
many of those referenced in the opening paragraph.

For this problem the weight function is $w(z) = 1/(2\pi i)$ and
the Cauchy transform is
\begin{equation}
\C(s) = \begin{cases} 1 & s \hbox{ interior to } \gamma, \\
0 & s \hbox{ exterior to } \gamma. \end{cases}
\label{cauchyjump}
\end{equation}
Thus $\C(s)$ has a jump of the simplest kind across~$\gamma$.
Let us suppose that $f$ is analytic in the annulus
$\Omega = \{z\in\complex:r^{-1} \le|z| \le r\}$ for some $r>1$.
(The results will depend little on the choice of $r$.)
As our enclosing contour $\Gamma$ of the kind shown on the right
in Figure~\ref{schematic}, we take
\begin{equation}
\Gamma = \riS \cup rS 
\end{equation}
with positive orientation on $rS$ and negative on $\riS$.
We now find a degree $n$ rational function $r_n$ that approximates 
$\C$ on $\Gamma$ as in (\ref{approxeps}).
The following code requests an approximation with accuracy $10^{-8}$:
\medskip

{\small
\begin{verbatim}
        S = exp(2i*pi*(1:100)'/100);
        Z = [2*S; 0.5*S];
        F = [zeros(size(S)); -ones(size(S))];
        [r,pol,res] = aaa(F,Z,'tol',1e-8,'sign',1,'lawson',20);
\end{verbatim}
\par}
\medskip

\begin{figure}
\begin{center}
\includegraphics[trim=45 94 20 7, clip, scale=.88]{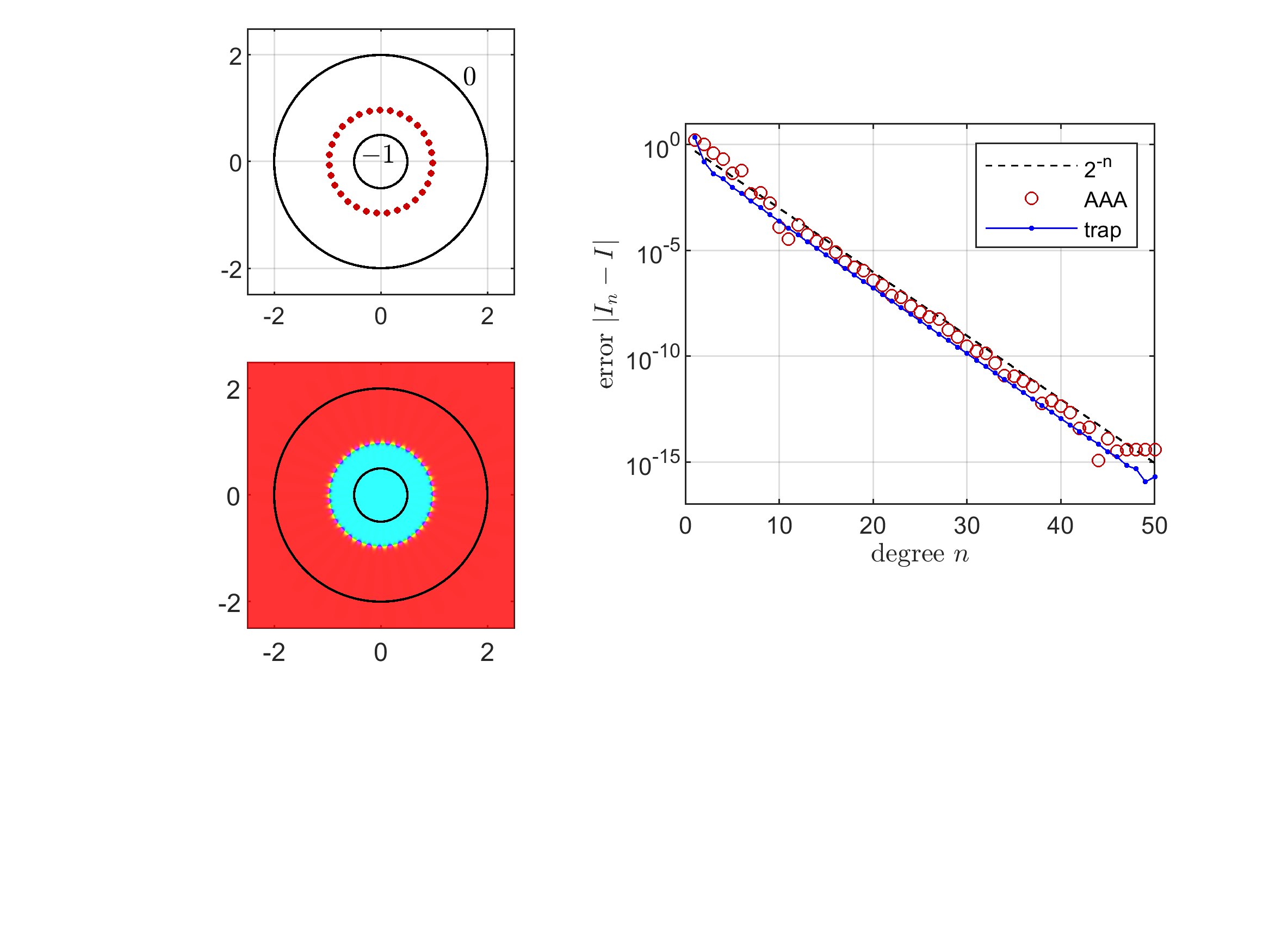}
\end{center}
\vspace*{-5pt}
\caption{\label{circlefig}On the left,
quadrature nodes---poles of a degree $n=31$ rational function---computed
in about $0.05$ s on a laptop by AAA rational approximation of
a two-branch target function with accuracy
specification $10^{-8}$.  The resulting approximate
integral\/ $(\ref{circint})$ 
has accuracy $2\cdot 10^{-10}$, as compared with
$7\cdot 10^{-11}$ for the exact\/ $31$-point trapezoidal rule.
The phase portrait now shows $r(z) + \half$ to separate values
$\approx -1$ and $\approx 0$.
On the right, errors as a function of\/ $n$.}
\end{figure}

\noindent The near-best approximation obtained has degree $n=31$,
with 31 poles ranging in modulus from $0.968$ to $0.970$, plotted
on the left in Figure~\ref{circlefig}.  The phase portrait of
$r(z) + \half$, with the constant $\half$ added to ensure a
separation between values $r(z) \approx -1$ and $r(z) \approx 0$,
shows bright cyan in the interior of the disk and bright red in
the exterior.  To check the accuracy of the approximation as a
quadrature formula, we apply it to (\ref{circint}) with integrand
$f(z) = -2\sqrt{z^2-1/4}$, which has a pair of branch points at
$z=\pm 1/2$.  The computed result {\tt In} is $0.999999999901 +
0.000000000144\kern .5pt i$, matching the correct value $1$ with an
error of about $2\cdot 10^{-10}$.  The exact 31-point trapezoidal
rule of (\ref{circint}) would improve this just a little bit to
$7\cdot 10^{-11}$.

\begin{figure}
\begin{center}
\includegraphics[trim=45 94 20 7, clip, scale=.88]{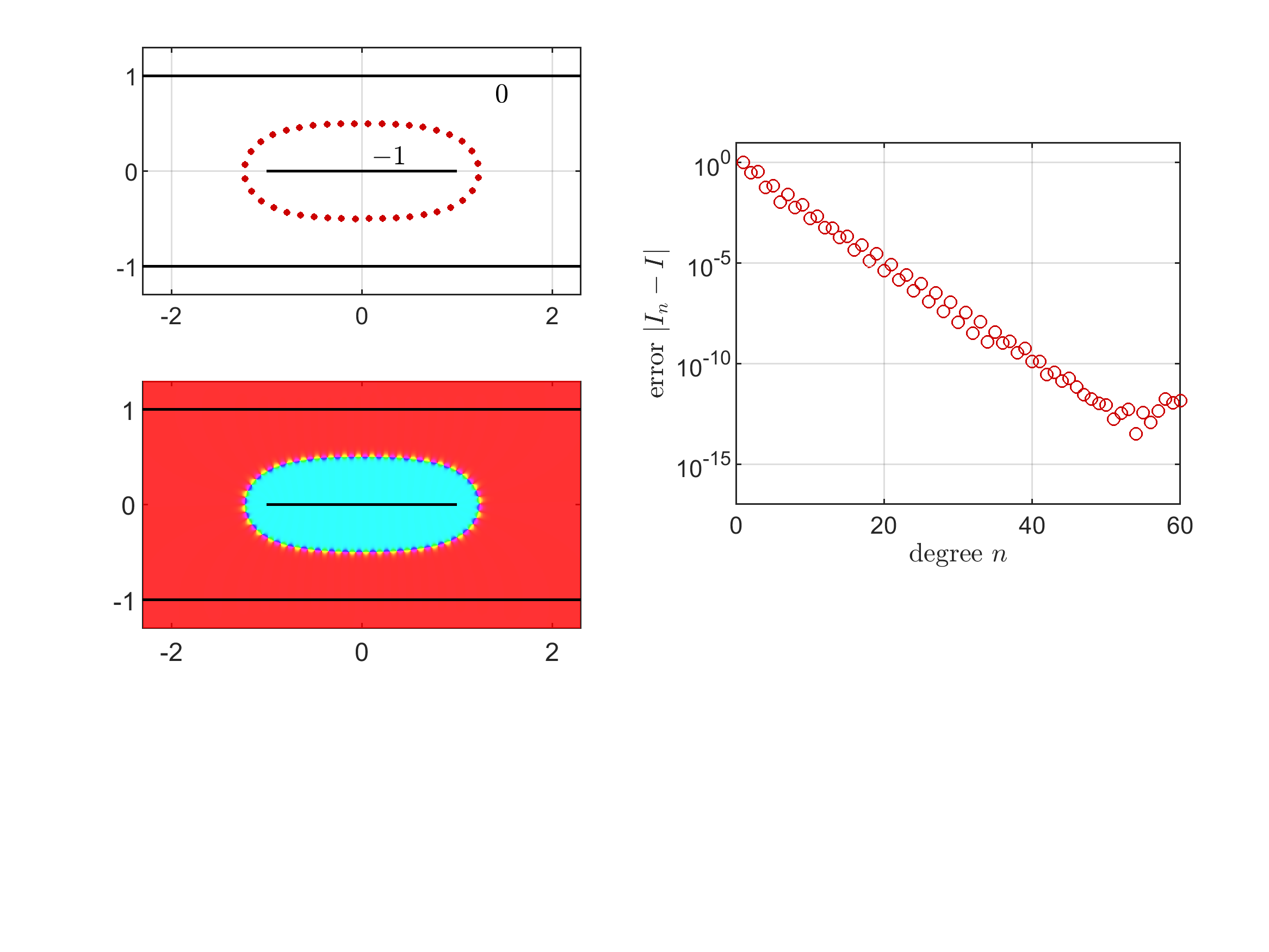}
\end{center}
\kern -5pt
\caption{\label{circlevfig}Repetition of Figure\/ $\ref{circlefig}$
based on the assumption that
$f$ is analytic in the strip 
$-1\le \Im z \le 1$ minus the unit interval $[-1,1]$.}
\end{figure}

For the variant problem of this section, shown
in Figure~\ref{circlevfig}, we suppose that it is known
that $f$ is analytic not in a circular annulus but in the region of
the strip $-1\le \Im z \le 1$ exterior to the unit interval $[-1,1]$.
One could derive a quadrature formula for
this geometry analytically by conformal mapping, but this would take some work,
and any modification of the geometry would require a new analysis.
Instead we proceed numerically and change the discretization above to
\medskip

{\small
\begin{verbatim}
        segment = linspace(1,1,200)';
        long = tan(pi*(-99:99)'/200);
        Z = [long+1i; long-1i; segment];
        F = [zeros(size([long; long])); -ones(size(segment))];
        [r,pol,res] = aaa(F,Z,'tol',1e-8,'sign',1,'lawson',20);
\end{verbatim}
\par}
\medskip

\noindent
The approximation comes out with degree 40, and we test it
on another function with integral $1$, $f(z) = -\sqrt{(z-1)/(z+1)}$.
The result is $0.999999999956 - 0.000000000024\kern .3pt
i$ with error $5.0\cdot 10^{-11}$.  

\section{Application 7: Functions of matrices via Cauchy integrals}

\begin{figure}
\begin{center}
\vspace{10pt}
\includegraphics[trim=0 0 0 0, clip,scale=0.21]{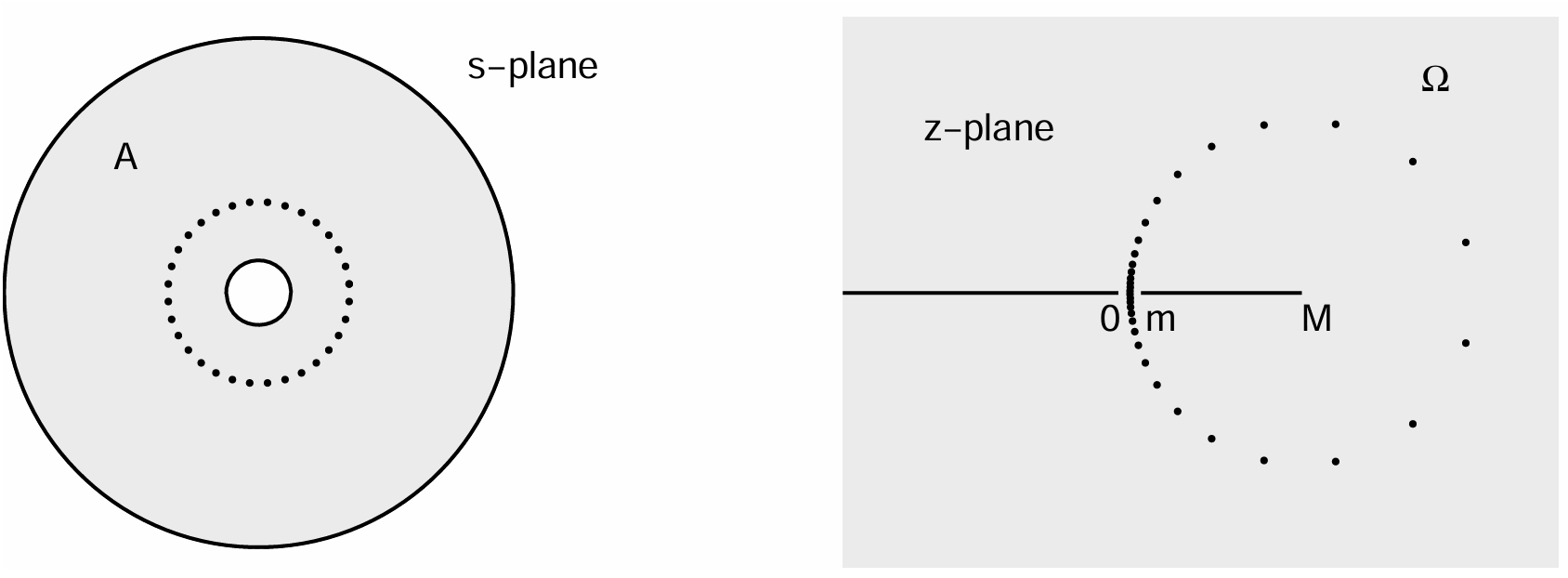}
\end{center}
\caption{\label{nicksfig}A figure from\/ {\rm \cite{nicks}} showing
conformal transplantation of the $32$-point equispaced trapezoidal rule in a circular annulus
to a domain
$\complex\kern 1pt\backslash \{ [m,M]\cup \negreal\}$ with
$m=1/8$, $M=1$.  This leads
to fast evaluation of functions like $A^\alpha$ or $\log(A)$ if $A$ is symmetric positive
definite with spectrum contained in $[m,M]$.}
\end{figure}

A function of a matrix
or operator $A$ can be defined by a Cauchy integral \cite{higham,kato},
\begin{equation}
f(A) = {1\over 2\pi i} \int_\Gamma f(z) (zI-A)^{-1} dz,
\label{dunford}
\end{equation}
where $\Gamma$ is a closed contour enclosing the spectrum of $A$
and contained in the region of analyticity of $f$.  This suggests
the use of discretizations to approximate $f(A) b$ for a given
vector $b$ by a linear combination of vectors $(\zk I-A)^{-1}
b$, where $\{\zk\}$ are quadrature nodes.  For example, suppose
$A$ is a symmetric positive definite matrix with eigenvalues in
$[m,M]$ for some $0<m<M<\infty$, as suggested in the right image of
Figure~\ref{nicksfig}, and suppose $f$ is a function like $A^\alpha$
or $\log(A)$ that is analytic in $\complex\kern 1pt\backslash
\negreal$.  Then the trapezoidal rule with $n$ equally spaced points
could be applied over a circle passing through the point $z=m/2$,
but the convergence as $n\to\infty$ would be very slow, requiring
$O(M/m)$ linear solves for each reduction of the error by a constant
factor.  In the ``three Nicks paper'' \cite{nicks}, it was shown that
this number can be improved enormously to $O(\kern .4pt \log(M/m))$
by a conformal map involving Jacobi elliptic functions, as shown in
the figure.\footnote{Zolotarev was an expert in elliptic functions,
and the mathematics here is very close to some of his work, though
so far as we are aware, he did not make connections between rational
approximation and quadrature.}  Now the mathematics is equivalent
to the application of the trapezoidal rule on a circular annulus
that is conformally equivalent to $\complex\backslash \{ [m,M]\cup
\negreal \}$.  This corresponds to a fast-converging transplanted
quadrature formula in unevenly spaced points in the $z$-plane.

To treat the problem by rational approximation instead
of conformal mapping, we note
that this is another case where the contour $\gamma$ is
flexible, to be determined, and the Cauchy transform will
again be (\ref{cauchyjump}).
Thus we can adjust the code segment of the last section to
\medskip

{\small
\begin{verbatim}
        segment = logspace(log10(1/8),0,100)';
        negreal = 1 - 1./linspace(.005,1,100)';
        Z = [negreal; segment];
        F = [zeros(size(negreal)); -ones(size(segment))];
        [r,pol,res] = aaa(F,Z,'degree',32,'sign',1,'lawson',0);
\end{verbatim}
\par}
\medskip

\begin{figure}
\begin{center}
\includegraphics[trim=45 94 20 7, clip, scale=.88]{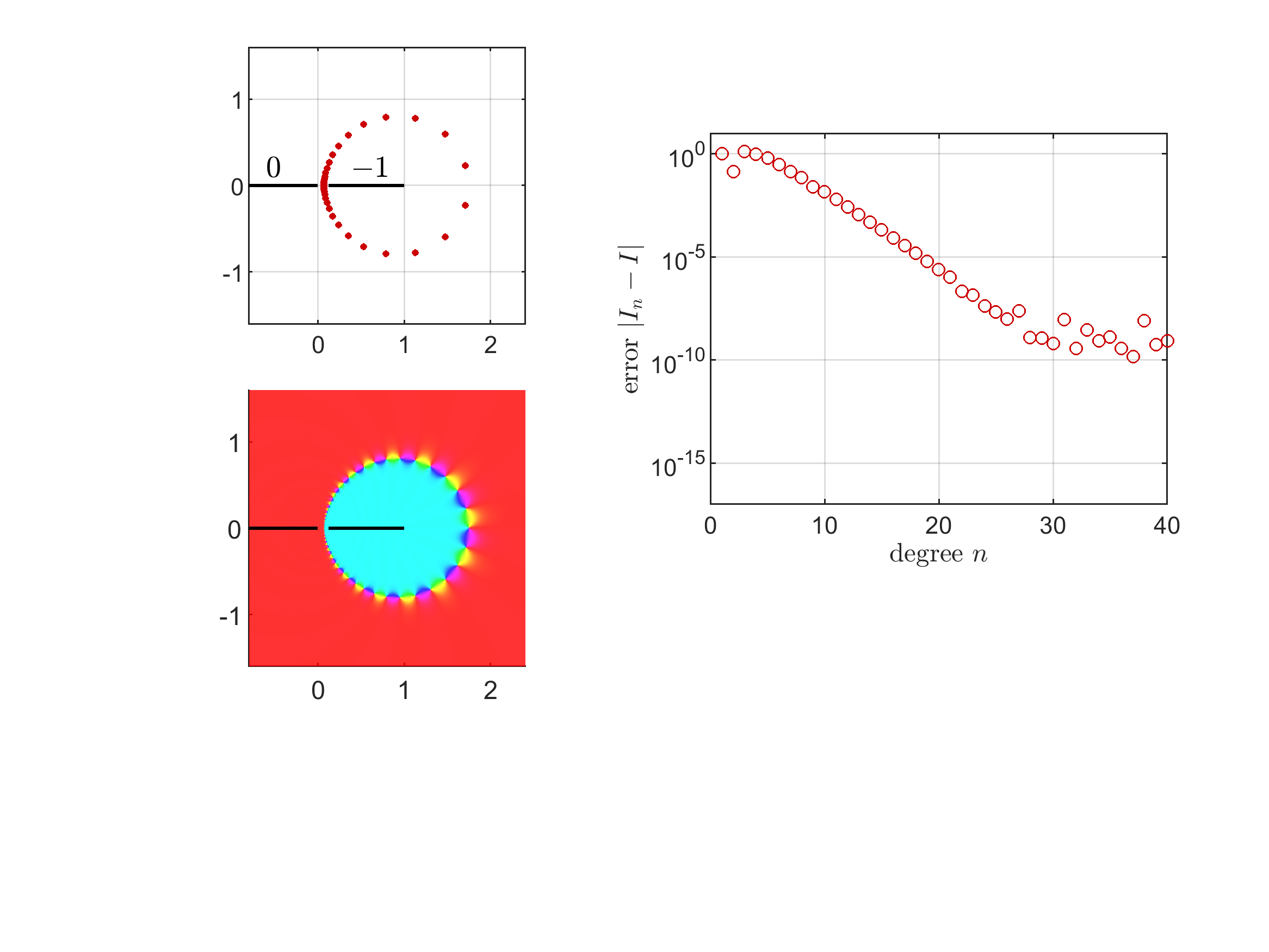}
\end{center}
\caption{\label{funAfig}Treatment of the same problem by rational approximation.  On
the left, quadrature poles with $n=32$, giving $8$-digit accuracy for a test
problem.  On the right, convergence for a test problem as a function of $n$.}
\end{figure}

\noindent  The poles shown on the left in Figure~\ref{funAfig}
are strikingly similar to those of Figure~\ref{nicksfig}, even though
neither a contour nor a node distribution has been fixed a priori.
Choosing another integrand with integral exactly 1,
$f(z) = (16/7)\sqrt{(z-1/8)/(z-1)}$,
gives the result $1.000000000092 - 0.000000000000\kern .4pt i$
with error magnitude $9.2\cdot 10^{-11}$.
This problem is more delicate than the last, and for clean
results as a function of degree $n$, shown on the right
in Figure~\ref{funAfig}, we modified the call to {\tt aaa.m} to
\medskip

{\small
\begin{verbatim}
     [r,pol,res] = aaa(F,Z,'degree',n,'sign',1,'lawson',50,'damping',.5);
\end{verbatim}
\par}
\medskip

\noindent
Such adjustments reflect the fact that with the rational
approximation algorithms of 2025, hand-tuning is often needed to get the
best results.  We hope that the situation will improve in the next
few years as the algorithms mature.

\begin{figure}
\begin{center}
\includegraphics[trim=45 195 36 0, clip, scale=.92]{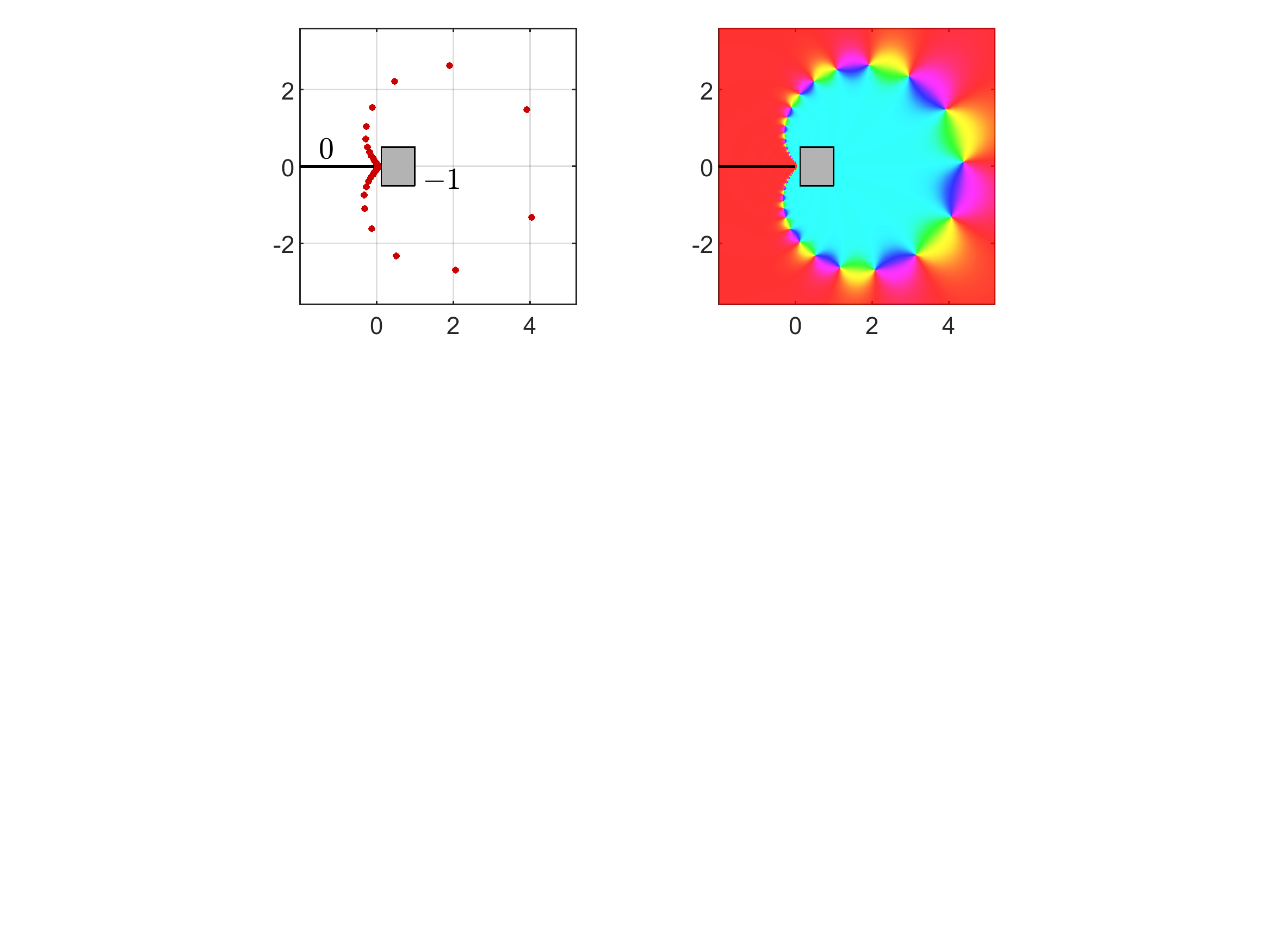}
\end{center}
\vspace*{-12pt}
\caption{\label{funAvfig}A variant of Figure $\ref{funAfig}$ in which
the spectrum of $A$ lies in a rectangle.}
\end{figure}

The variant problem of this section was suggested to us by Maria
L\'opez-Fern\'andez of the University of Malaga.  What if 
$A$ is not symmetric but is known to have its spectrum
in a given region of the right half-plane, like a rectangle?
Figure~\ref{funAvfig} shows that the approximation $r$ now has
poles extending much further away from the origin, forming a bend
of close to $90^\circ$ near the spectrum.

\section{Discussion}

Nothing is entirely new in the field of quadrature, and the link
with rational functions starts with Gauss in 1814 \cite{gauss}
and has been exploited in many ways since then.
Most of this work has been algebraic, however, related to Pad\'e
approximations at $\infty$ and associated orthogonal polynomials.
This paper has followed a more analytic approach based on numerical
rational approximations made possible by the AAA algorithm.
Here there is a smaller tradition, including the contributions
\cite{tm,gausscc,talbot}.

Whenever you see a string of quadrature nodes, you should consider
it as a rational approximation to a branch cut.  That is our
conceptual message.  Our practical message is that it is amazingly
easy to compute near-optimal quadrature formulas by taking advantage
of this connection.  We make no claim that any of the particular
techniques we have illustrated will immediately bring improvements
in computational science, but we have outlined many possibilities
and are confident that some of them will prove useful.

Jumps of complex functions across arcs arise in many contexts
and are associated with notions including Cauchy/Stieltjes/Hilbert
transforms, the Sokhotski--Plemelj formula, Wiener--Hopf factorization,
the Dirichlet-to-Neumann map, and Rie\-mann--Hilbert problems.
The literature of such topics is enormous, and links could be found
from this paper to most of it.  Another topic in this space
is that of {\em hyperfunctions,} which are generalized functions
realized as (equivalence classes of) pairs of analytic functions
on either side of an arc \cite{graf}.  The present paper could be
said to investigate quadrature formulas as rational hyperfunction
approximations to Cauchy transforms.   Perhaps this work can help
encourage the growth of a numerical theory of hyperfunctions,
which has been mostly a theoretical topic heretofore.

Theoretical questions lie around every corner.  Exactly what are the
links with potential theory?  Does the phenomenon of ``balayage'' \cite{gust}
relate interpolation points at~$\infty$, for example, to almost
equally effective points along a Bernstein ellipse?  Can backward
error analysis confirm our experience that nodes and weights may
be individually far from optimal yet still combine to give nearly
optimal quadratures?  
Can one prove in some generality that best rational approximations
match the convergence rates of Gauss quadrature rules?
When can weights (= residues) be guaranteed
to be positive, or at least uniformly summable as $n\to\infty$?
These are just a few examples of many problems to be explored.

\section*{Acknowledgments}
We are grateful for helpful suggestions from Matt Colbrook,
Toby Driscoll, Stefan G\"uttel, Nick Hale, Arno Kuijlaars, 
Maria L\'opez-Fern\'an\-dez, Yuji Nakatsukasa, 
Sheehan Olver, Mika\"el Slevinsky, Tom Trogdon, Kieran Vlahakis, Andr\'e 
Weideman, and Heather Wilber.

\end{document}